\documentclass[english,reqno]{amsart}
\usepackage[T1]{fontenc}
\usepackage[latin9]{inputenc}
\usepackage{color}
\usepackage{babel}
\usepackage{verbatim}
\usepackage{amsthm}
\usepackage{amsmath}
\usepackage{amssymb}
\PassOptionsToPackage{normalem}{ulem}
\usepackage{ulem}
\usepackage[unicode=true]
 {hyperref}
\usepackage{breakurl}

\makeatletter
\numberwithin{equation}{section}
\numberwithin{figure}{section}
\theoremstyle{plain}
\newtheorem{thm}{\protect\theoremname}[section]
  \theoremstyle{definition}
  \newtheorem{defn}[thm]{\protect\definitionname}
  \theoremstyle{remark}
  \newtheorem{rem}[thm]{\protect\remarkname}
  \theoremstyle{plain}
  \newtheorem{cor}[thm]{\protect\corollaryname}
  \theoremstyle{plain}
  \newtheorem{fact}[thm]{\protect\factname}
  \theoremstyle{plain}
  \newtheorem{prop}[thm]{\protect\propositionname}
  \theoremstyle{definition}
  \newtheorem{example}[thm]{\protect\examplename}
  \theoremstyle{definition}
  \newtheorem{xca}[thm]{\protect\exercisename}

\usepackage{amscd}

\theoremstyle{plain}

\makeatother

  \providecommand{\corollaryname}{Corollary}
  \providecommand{\definitionname}{Definition}
  \providecommand{\examplename}{Example}
  \providecommand{\exercisename}{Exercise}
  \providecommand{\factname}{Fact}
  \providecommand{\propositionname}{Proposition}
  \providecommand{\remarkname}{Remark}
\providecommand{\theoremname}{Theorem}

\begin{document}

\title[Riemann Integration]{Lecture notes about a simpler approach to Riemann integration}

\author{Michael Cwikel}

\address{Technion - Israel Institute of Technology, 32000 Haifa, Israel}

\email{mcwikel@math.technion.ac.il}

\maketitle

\section{A preface to students and fellow teachers}

\textcolor{red}{}%

This document will surely need further revisions and additions, and
some (hopefully only a few) corrections. But I consider that this
preliminary version can already be useful to some of you.

I hereby invite you to consider teaching and learning about the Riemann
integral of functions of one and more variables in a new way!~\textcolor{red}{}%
In fact it is not really new. It is no secret (see e.g., \cite{ST}
and \cite{R} pp.~50-51) that there is a considerably simpler and
more elegant way of defining and studying these integrals, than the
way that most people and books have been using for many years. But
somehow no one, as far as I have been able to discover, has actually
bothered to use this way systematically in some course or textbook.
\textcolor{red}{}%

Please give me a few moments to indicate an equivalent simpler definition
of the Riemann integral $\int_{a}^{b}f(x)dx$. We will not need to
mess around here with an (uncountable!!) family of different partitions
of the interval $[a,b]$. \textcolor{red}{}%
{} We will only have to write down two \textit{explicitly defined} sequences
of numbers $\{U_{k}\}$ and $\left\{ L_{k}\right\} $. It is easy
to prove \textcolor{red}{}%
{} that they both converge. If they converge to the same limit then
the above integral exists and equals that limit. \textcolor{red}{}%

In any course where it is necessary to prove the various basic properties
of this integral and of integrable functions (integrability of their
sums and products and of $\max\left\{ f,g\right\} $ etc., the interaction
with continuity) then most of the proofs will be somewhat simpler
than the usual ones. 

But when one comes to deal with double integrals and triple integrals
there are even greater advantages to be harvested. (It's a pity that
the authors of \cite{ST} did not take time to point this out.)

If we use appropriate notation, the definition of the $n$-fold integral
of a function $f$ of $n$-variables on some subset of $\mathbb{R}^{n}$
will be \textit{exactly} (or almost exactly) the same as the one given
above for $n=1$. And the proofs of the properties of these integrals
will be almost exactly the same. 

In some versions of the ``classical'' approach one first defines
the double integral only on a rectangle, considering all possible
decompositions of that rectangle into smaller (``parallel'') rectangles.
Then one has to talk about more general sets which have ``area''
(are Jordan measurable) and describe how to define the double integral
of a function on such sets. One needs to impose some conditions on
the function and some other conditions on the set. By the time one
gets to triple integrals, both the teacher and the students have often
run out of energy and patience and time, and content themselves with
generally claiming that the theory for triple integrals is essentially
the same. (I confess to making some analogous claims here, but I think
that the exercises needed to justify them, which can use use tools
prepared here, are rather shorter and simpler than in some classical
treatments.)

In the approach that I suggest we should adopt, we bypass a lot of
this messing around. If an engineer has asked us to calculate or think
about the $n$-fold integral of some function $f$ on some subset
$E$ of $\mathbb{R}^{n}$ then we simply consider the function $g$
on $\mathbb{R}^{n}$ which equals $f$ on $E$ and $0$ everywhere
else, and apply exactly the definition via the two explicit sequences
mentioned above (and formulated below) to the new function $g$.

Of course we also want to explicitly calculate $n$-fold integrals
by reducing them to repeated integrals (a ``baby'' version of Fubini's
theorem) and here also (see Theorem \ref{thm:fubini1-1} and Exercise
\ref{ex:Fubini} below) I can offer you a rather simple formulation
with a simple proof of this baby Fubini theorem. It will readily apply
to the usual exercises in calculus courses which deal with double
or triple integrals. 

I have not sought to consider ALL Riemann integrable functions in
presenting the material below, since I did not wish to deal with (outer)
Lebesgue measure. However functions of compact support which are continuous
(and therefore uniformly continuous) are of course shown to be integrable,
and it is not too too complicated to extend that result to show that
that any bounded function of compact support whose points of discontinuity
are contained in a ``very small'' set is also Riemann integrable.
For $E$ to be ``very small'' means that $E$ has Jordan measure
$0$, but this can be equivalently defined more comfortably by the
condition that the explicit sequence $\left\{ U_{k}\right\} $ for
the characteristic function of $E$ converges to $0$.

So for most of our purposes here we do not need to define Jordan measure
in general. But anyone who needs to know about Jordan measure and
Jordan measurability of other bounded subsets of $\mathbb{R}^{n}$
can get what they need by simply considering the integral of the characteristic
functions of such sets. This will drop out easily as a by-product
of the theory, instead of being a messy pre-requisite for developing
the theory.

If we want to reap the full advantages of using this approach, we
should, if possible, delay treatment of all integrals until the students
have at least a minimal knowledge of functions of two and three variables,
including the notion of continuity. So this could require juggling
the order of some topics within and between standard existing calculus
or analysis courses. Such a change does not render existing fine textbooks
and lecture notes obsolete. We simply ask students to read most sections
of them in a different order. And a small number of sections from
those sources should be replaced by new material, along the lines
of what appears here in the following sections.

If it turns out that other courses (such as physics) that some students
are taking in parallel with our calculus or analysis course require
some earlier familiarity with integrals, we can insert a lecture or
two early in our course where we intuitively define and intuitively
explain some physical and geometric interpetations of the definite
integral $\int_{a}^{b}f(x)dx$ only in the special case where $f$
is continuous. In that lecture we can state and ask the students to
temporarily believe a small number of the properties of such integrals,
which of course we will be obliged to prove later in the course when
we extend the definition of integrals to a more general situation.
These properties should suffice for early physics needs, and also
for any results about (uniformly convergent) sequences of functions
which arise in our course. They can also be used to give a rigorous
proof (or explanation) that $\frac{d}{dx}\int_{a}^{x}f(t)dt=f(x)$
at each point of continuity $x$ of an integrable function $f$ and
from this one can readily deduce that $\int_{a}^{b}g'(t)dt=g(b)-g(a)$
for every function $g$ which has a continuous derivative in $\left[a,b\right]$.
(Cf.~Subsection \ref{sub:meo} below.) \textcolor{red}{}%

Of course it is natural to ask whether a student who is destined in
any case to ultimately study the Lebesgue integral should be required
to study the Riemann integral in detail. I would suggest that even
if some mathematics departments decide to do nothing more than give
their students a brief survey of Riemann integration, then they can
do so more easily and effectively by using the definitions which appear
in this document. Of course the classical definitions remain important
for those experts who need them as a first step towards the intricacies
of the Henstock-Kurzweil integral etc.

Here are the main aims of the following material. Some of them, in
particular the seventh, remain to be dealt with in future versions
of this document. 

1. To define integrable functions of one variable and their integrals
via dyadic intervals. Then to show how to extend these definitions
to functions of 2 or 3 or $m$ variables, via dyadic squares or cubes.

2. To obtain some basic properties of integrals and integrable functions.

3. To show how to get new integrable functions by combining other
functions in various ways, when these other functions are already
known to be integrable.

4. To show that bounded functions on a bounded set which are continuous
everywhere in that set, except perhaps on some ``very small'' subset,
are integrable functions.

5. To show how to calculate some integrals of functions of one variable
via primitive functions.

6. To show how to calculate some integrals of functions of $m$ variables
via ``repeated integration''.

7. To introduce and study reasonable definitions of length, area and
volume of sufficiently ``nice'' subsets of $\mathbb{R}$, $\mathbb{R}^{2}$
and $\mathbb{R}^{3}$ respectively. 

8. To show that the definitions of integrability and integrals given
here via dyadic intervals or squares or cubes are equivalent to the
``classical'' definitions via arbitrary partitions.

\medskip{}

The presentation here moves at a gentle pace. Apologies to those of
you who find some of the explanations superfluous.

\textit{If you happen to have seen a book or set of lecture notes
which uses this approach, please please tell me about it. It would
be quite surprising if such a formal document does not yet exist somewhere.}

\section{Notation, terminology and recalling of some standard results.}

A left semiclosed interval is an interval of the form $[a,b)$ for
$-\infty<a<b<\infty$.

For any bounded interval $I\subset\mathbb{R}$ we let $\left|I\right|$
or $\left|I\right|_{1}$ denote the length of $I$. Here $I$ may
be open, closed or semiclosed. 

For any $m\in\mathbb{N}$ and for any subset $E$ of $\mathbb{R}^{m}$
we denote the interior, closure and boundary of $E$ respectively
by $E^{\circ}$, $\overline{E}$ and $\partial E$.

For any subset $H$ of $\mathbb{R}^{m}$ and any function $f:H\to\mathbb{R}$,
we recall that $f$ is \textbf{\textit{uniformly continuous on }}$H$
if, for each $\epsilon>0$, there exists a number $\delta(\epsilon)>0$
such that, whenever two points $p=\left(p_{1},p_{2},...,p_{m}\right)$
and $q=\left(q_{1},q_{2},....,q_{m}\right)$ in $H$ satisfy $\max\left\{ \left|p_{j}-q_{j}\right|:j=1,2,...,m\right\} <\delta(\epsilon)$,
then $\left|f(p)-f(q)\right|<\epsilon$. Of course in this definition
you can, if you wish, equivalently replace $\max\left\{ \left|p_{j}-q_{j}\right|:j=1,2,...,m\right\} $
(which is often denoted by $d_{\infty}(p,q)$) by the usual euclidean
distance $d_{2}(p,q)=\sqrt{\sum_{j=1}^{m}(p_{j}-q_{j})^{2}}$, since
$d_{\infty}(p,q)\le d_{2}(p,q)\le\sqrt{m}d_{\infty}(p,q)$. 

We recall the result that if $H$ is a bounded and closed subset of
$\mathbb{R}^{m}$, then any function $f:H\to\mathbb{R}$ which is
continuous at every point of $H$ is uniformly continuous on $H$.
(Here continuity at a point means continuity only with respect to
$H$. This means that if, for example $m=1$ and $H$ is a closed
interval, then $f$ has to be one-sidedly continuous at the endpoints
of $H$.)

Suppose that $m\in\mathbb{N}$ and $E_{1}$, $E_{2}$, ...., $E_{m}$
are all subsets of $\mathbb{R}$. Then their \textbf{\textit{cartesian
product}} is the set 
\begin{equation}
E_{1}\times E_{2}\times...\times E_{m}:=\left\{ \left(x_{1},x_{2},...,x_{m}\right)\in\mathbb{R}^{m}:x_{j}\in E_{j},\,\,1\le j\le m\right\} .\label{eq:dcp}
\end{equation}
We will also use the notation $\prod_{j=1}^{m}E_{j}$ for this set. 

For each $m\in\mathbb{N}$ we define an \textbf{\textit{$m$-rectangle}}
to be a set of the form $E=\prod_{j=1}^{m}I_{j}$ where each $I_{j}$
is a bounded interval in $\mathbb{R}$. (Each $I_{j}$ can be open,
closed or semiclosed.) If all the intervals $I_{j}$ have the same
length we can also refer to $E$ as an $m$\textbf{\textit{-cube}},
or simply a cube, if the value of $m$ is clear from the context.
We define the \textbf{\textit{$m$-dimensional volume}} of the above
set $m$-rectangle $E$ to be the number $\left|E\right|_{m}=\prod_{j=1}^{m}\left|I_{j}\right|_{1}$.
This means that $1$-rectangles are simply intervals in $\mathbb{R}$,
$2$-rectangles are simply rectangles in $\mathbb{R}^{2}$ and $3$-rectangles
are simply ``boxes'' (sometimes called rectangular cuboids or rectangular
parallelepipeds) in $\mathbb{R}^{3}$. In each case, of course, the
sides %
{} of these sets are parallel to the coordinate axes. For $m=1,2,3$
the $m$-dimensional volume of an $m$-rectangle is, respectively,
its length, its area, or its volume in the usual sense of the word.

When $m>3$ it is quite difficult to find a way to draw ``pictures''
of $m$-rectangles or of any other subsets of $\mathbb{R}^{m}$. But
the theory that we want to describe can be developed just as easily,
or almost as easily for every choice of $m$, as for the easier to
visualize cases where $m$ is $1$ or $2$ or $3$.

For each $m\in\mathbb{N}$ and each set $E\subset\mathbb{R}^{m}$,
we use the notation $\chi_{E}$ to denote the \textbf{\textit{characteristic
function}}, or \textbf{\textit{indicator function of}} $E$. This
is a real valued function defined on all of $\mathbb{R}^{m}$ by 
\[
\chi_{E}(p)=\left\{ \begin{array}{ccc}
1 & , & p\in E\\
\\
0 & , & p\in\mathbb{R}^{m}\setminus E
\end{array}\right..
\]
(In some documents the notation $\mathbf{I}_{E}$ is used instead
of $\chi_{E}$.) 

Suppose that $E$ and $G$ are subsets of $\mathbb{R}^{m}$ with $E\subset G$
and we are given a function $f:G\to\mathbb{R}$. Let us agree to adopt
the convention that $f\chi_{E}$, the product of the functions $f$
and $\chi_{E}$, will in fact denote a function defined on all of
$\mathbb{R}^{m}$, not just on $G$, which is defined to equal $f$
at all points of $E$ and to equal $0$ at all points of $\mathbb{R}^{m}\setminus E$.
This could be considered an ``abuse of notation'', since $f$ is
possibly not defined at some points of $\mathbb{R}^{m}\setminus E$.
But it will cause no difficulties or ambiguities.

\section{\label{sec:SupAndInf}Some basic properties of suprema and infima}

Let $A$ be a bounded set of real numbers. In case anyone has forgotten,
the numbers $\sup A$ and $\inf A$ are the (necessarily unique) numbers
$a$ and $b$ having the properties that $x\le a$ and $x\ge b$ for
every $x\in A$ and that, for each $n\in\mathbb{N}$ there exist points
$a_{n}$ and $b_{n}$ in $A$ such that $a-a_{n}\le1/n$ and $b_{n}-b\le1/n$.

Suppose that $E$ is a subset of $\mathbb{R}$ or $\mathbb{R}^{n}$
and that $f$ and $g$ are real valued function defined on $E$. Then
it is easy to check the following facts which will be used in the
discussion to follow:

\begin{equation}
\mbox{ }\begin{array}{ccc}
 & {\displaystyle \sup\left\{ cf(p):p\in E\right\} =c\sup\left\{ f(p):p\in E\right\} } & \phantom{\begin{array}{c}
x\\
.
\end{array}}\\
\mbox{and} & {\displaystyle \inf\left\{ cf(p):p\in E\right\} =c\inf\left\{ f(p):p\in E\right\} } & \phantom{\begin{array}{c}
.\\
x
\end{array}}
\end{array}\mbox{ for every constant }c\ge0.\label{eq:sp01}
\end{equation}

\smallskip{}

\begin{equation}
\begin{array}{ccc}
 & {\displaystyle \sup\left\{ cf(p):p\in E\right\} =c\inf\left\{ f(p):p\in E\right\} } & \phantom{\begin{array}{c}
x\\
.
\end{array}}\\
\mbox{and} & {\displaystyle \inf\left\{ cf(p):p\in E\right\} =c\sup\left\{ f(p):p\in E\right\} } & \phantom{\begin{array}{c}
.\\
x
\end{array}}
\end{array}\mbox{ for every constant }c\le0.\label{eq:sp02}
\end{equation}

\smallskip{}

\begin{equation}
\sup\left\{ f(p)+g(p):p\in E\right\} \le\sup\left\{ f(p):p\in E\right\} +\sup\left\{ g(p):p\in E\right\} .\label{eq:sp-add}
\end{equation}

\begin{equation}
\inf\left\{ f(p)+g(p):p\in E\right\} \ge\inf\left\{ f(p):p\in E\right\} +\inf\left\{ g(p):p\in E\right\} .\label{eq:sp-add-inf}
\end{equation}

For every subset $G\subset E$,
\begin{equation}
\sup\left\{ f(p):p\in G\right\} \le\sup\left\{ f(p):p\in E\right\} \mbox{ and }\inf\left\{ f(p):p\in G\right\} \ge\inf\left\{ f(p):p\in E\right\} `.\label{eq:sp-inc}
\end{equation}

If $f(p)\le g(p)$ for every $p\in E$ then 
\begin{equation}
\sup\left\{ f(p):p\in E\right\} \le\sup\left\{ g(p):p\in E\right\} \mbox{ and }\inf\left\{ f(p):p\in E\right\} \le\inf\left\{ g(p):p\in E\right\} .\label{eq:sp-mon}
\end{equation}

We will usually use the notation $\sup_{p\in E}f(p)$ instead of $\sup\left\{ f(p):p\in E\right\} $
and $\inf_{p\in E}f(p)$ instead of $\inf\left\{ f(p):p\in E\right\} $.\textit{ }

We will also need to know that 
\begin{equation}
\sup_{p\in E}f(p)-\inf_{p\in E}f(p)=\sup\left\{ \left|f(p)-f(q)\right|:p,q\in E\right\} \label{eq:burz}
\end{equation}
where the right hand side of (\ref{eq:burz}) can also be denoted
by $\sup_{p,q\in E}\left|f(p)-f(q)\right|$. 

The proof of (\ref{eq:burz}) is quite straightforward, but let us
present it anyway: We first note that, for each $p$ and $q$ in $E$,
\[
\left|f(p)-f(q)\right|=\max\left\{ f(p),f(q)\right\} -\min\left\{ f(p),f(q)\right\} \le\sup_{p\in E}f(p)-\inf_{p\in E}f(p).
\]
Taking the supremum over all choices of $p$ and $q$ in $E$ gives
us ``half'' of (\ref{eq:burz}), i.e., the inequality ``$\ge$''.
In order to obtain the reverse inequality ``$\le$'', we let $\left\{ p_{n}\right\} _{n\in\mathbb{N}}$
and $\left\{ q_{n}\right\} _{n\in\mathbb{N}}$ be sequences of points
in $E$ such that $\sup_{p\in E}f(p)=\lim_{n\to\infty}f(p_{n})$ and
$\inf_{p\in E}f(p)=\lim_{n\to\infty}f(q_{n})$. Then 
\[
\sup_{p\in E}f(p)-\inf_{p\in E}f(p)=\lim_{n\to\infty}\left(f(p_{n})-f(q_{n})\right)\le\lim_{n\to\infty}\left|f(p_{n})-f(q_{n})\right|
\]
and this last limit does not exceed $\sup_{p,q\in E}\left|f(p)-f(q)\right|$.
This completes the proof of (\ref{eq:burz}).

\section{\label{sec:BasicProps}Definitions and basic properties of integral
and integrable functions}

We ultimately want to deal with functions of $m$ variables. 

Let us start with functions of one variable, i.e., $m=1$, but we
will be setting things up so that the transition to functions of $2$
or $3$ or even $n$ variables can be an almost triviality.

You have given me a bounded function $f$ defined on some bounded
subset $E$ of $\mathbb{R}$. You want me to define, and or calculate
the integral of $f$ on $E$. (Of course $E$ will usually be an interval.)

With your kind permission, I will think of $f$ as being defined on
all of $\mathbb{R}$ by setting it equal to $0$ outside $E$.

It will be convenient to let $BB(\mathbb{R})$ denote the set of all
functions $f:\mathbb{R}\to\mathbb{R}$ which are bounded and which
vanish outside some bounded set. 

Let $D_{0}$ be the set of all intervals $I$ of the special form
$I=[n,n+1)$ where $n\in\mathbb{Z}$.

More generally, for each non-negative integer $k$, let $D_{k}$ be
the set of all intervals $I$ of the special form $I=[n2^{-k},n2^{-k}+2^{-k})$
where $n\in\mathbb{Z}$.

Now, for each $f\in BB(\mathbb{R})$, and each $k\in\mathbb{N}\cup\left\{ 0\right\} $,
define $U_{k}(f)=2^{-k}\sum_{I\in D_{k}}\sup_{p\in I}f(p)$ and $L_{k}(f)=2^{-k}\sum_{I\in D_{k}}\inf_{p\in I}f(p)$. 

(Of course the letters $U$ and $L$ here stand for ``upper'' and
``lower''.)

Of course here it would be good to draw some pictures, and maybe even
use some computer animations.

\textit{In ``low power'' courses which present fewer proofs and
where the students do not feel comfortable with suprema and infima.
these can be replaced by maxima and minima of the functions (which
can be assumed to be continuous) on the corresponding closed intervals. }

\textit{If you are wondering why I prefer to use semiclosed intervals
instead of closed intervals (as is done in the version of this approach
in \cite{ST}) you can find one reason for this by looking below at
the proof of a ``baby Fubini theorem'' Theorem \ref{thm:fubini1-1}.
For anyone who needs to know that it is equivalent to work with closed
or semiclosed (or open) intervals, this can be proved fairly easily.
(Cf.~Section \ref{sec:Equivalence} for a proof of a more general
related result.) But we don't need to worry about that just now.}

Of course in each of these sums $U_{k}(f)$ and $L_{k}(f)$ there
will only be finitely many non-zero terms. Since each interval $I$
in $D_{k}$ is the union of two ( $=2^{m}$ ) adjacent intervals of
$D_{k+1}$ which we can denote by $J_{1}(I)$ and $J_{2}(I),$ it
is easy to see that 
\begin{equation}
2^{-(k+1)}\left(\sup_{p\in J_{1}(I)}f(p)+\sup_{p\in J_{2}(I)}f(p)\right)\le2^{-(k+1)}\left(\sup_{p\in I}f(p)+\sup_{p\in I}f(p)\right)=2^{-k}\sup_{p\in I}f(p).\label{eq:blrp}
\end{equation}
Summing this inequality over all $I\in D_{k}$ (so that $J_{1}(I)$
and $J_{2}(I)$ range over all intervals in $D_{k+1}$) gives us that
\begin{equation}
U_{k+1}(f)\le U_{k}(f).\label{eq:rwq}
\end{equation}
A similar argument shows that 
\begin{equation}
L_{k+1}(f)\ge L_{k}(f).\label{eq:eeq}
\end{equation}
(Perhaps you might prefer to deduce (\ref{eq:eeq}) from (\ref{eq:rwq})
and the rather obvious fact that $L_{k}(f)=-U_{k}(-f)$.) Also $L_{0}(f)\le L_{k}(f)\le U_{k}(f)\le U_{0}(f).$
So these two monotonic sequences $\left\{ U_{k}(f)\right\} _{k\in\mathbb{N}}$
and $\left\{ L_{k}(f)\right\} _{k\in\mathbb{N}}$ are bounded and
therefore converge. 
\begin{defn}
\label{def:dori}Let $f$ be a function in $BB(\mathbb{R})$. If $\lim_{k\to\infty}U_{k}(f)=\lim_{k\to\infty}L_{k}(f)$
then we say that $f$ is Riemann integrable on $\mathbb{R}$ and the
common value of these two limits is called the Riemann integral of
$f$ on $\mathbb{R}$ and may be denoted by $\int_{\mathbb{R}}f$.
\end{defn}
It will be convenient to let $RI(\mathbb{R})$ denote the set of all
Riemann integrable functions $f:\mathbb{R}\to\mathbb{R}$. There are
other more ``classical'' kinds of notation for $\int_{\mathbb{R}}f$.
In the more traditional treatment of Riemann integration (as for example
in \cite{R} or \cite{S}) we consider a function which is defined
only on some bounded interval $\left[a,b\right]$ and often use the
notation $\int_{a}^{b}f(x)dx$ for its integral on $\left[a,b\right]$.
It can be proved that this integral, however, is exactly the same
as the integral $\int_{\mathbb{R}}f$ defined just above, provided
we extend $f$ to all of $\mathbb{R}$ by setting $f(x)=0$ for all
$x\in\mathbb{R}\setminus\left[a,b\right]$. 

\textit{If you are worried that Definition \ref{def:dori} is not
equivalent to the usual (in my opinion rather messier) definition,
you can see a proof of this equivalence in Section \ref{sec:Equivalence}.
There is also a rather shorter ``almost'' proof of this equivalence
on the last page of \cite{ST}. I say ``almost'' because that proof
needs some adjustment because here we use semiclosed intervals instead
of closed intervals. It could perhaps be claimed that our students,
at least at this stage, simply do not need to relate to or know that
other messier definition. }

The intervals in the collections $D_{k}$ are usually called \textbf{\textit{dyadic
intervals}}. When we come to deal with double, or triple (or $n$-fold)
integrals, we simply replace these dyadic intervals by \textbf{\textit{dyadic
squares}}, or \textbf{\textit{dyadic cubes}}. So $D_{k}(\mathbb{R}^{2})$
is the collection of all squares of the form $[n_{1}2^{-k},n_{1}2^{-k}+2^{-k})\times[n_{2}2^{-k},n_{2}2^{-k}+2^{-k})$
and $D_{k}(\mathbb{R}^{3})$ is the collection of all cubes of the
form $[n_{1}2^{-k},n_{1}2^{-k}+2^{-k})\times[n_{2}2^{-k},n_{2}2^{-k}+2^{-k})\times[n_{3}2^{-k},n_{3}2^{-k}+2^{-k})$
where $n_{1},n_{2},n_{3}\in\mathbb{Z}.$

It is easy to draw a picture of $\mathbb{R}^{2}$ covered by \textcolor{red}{}%
``tiles'' of side length 1, (the squares of $D_{0}(\mathbb{R}^{2})$,
then show how each of these squares is split into four disjoint squares
of side length 1/2 in $D_{1}\left(\mathbb{R}^{2}\right)$, and so
on. 

The splitting of $\mathbb{R}^{3}$ into boxes/cubes in $D_{0}(\mathbb{R}^{3})$
and then the splitting of each such box into eight disjoint boxes
in $D_{1}\left(\mathbb{R}^{3}\right)$, and so on, is harder to draw,
but not harder to understand. 

So now if $m=1,2$ or $3$ (or even beyond) and if $f:\mathbb{R}^{m}\to\mathbb{R}$
is a bounded function of $m$ variables vanishing outside some bounded
set, (i.e., $f\in BB\left(\mathbb{R}^{m}\right)$) then we simply
set 

\begin{center}
$U_{k}(f)=2^{-k{\color{red}m}}\sum_{I\in D_{k}{\color{red}(\mathbb{R}^{m})}}\sup_{p\in I}f(p)$
and $L_{k}(f)=2^{-k{\color{red}m}}\sum_{I\in D_{k}{\color{red}(\mathbb{R}^{m})}}\inf_{p\in I}f(p)$. 
\par\end{center}

These are almost the same as the formulae given above for $m=1$.
(You can presumably see the small changes which I have written in
red.) Sometimes, when it is desirable to more explicitly indicate
that we are dealing with functions of $m$ variables, we may use the
notation $U_{k}(f,\mathbb{R}^{m})$ instead of $U_{k}(f)$ and similarly
$L_{k}(f,\mathbb{R}^{m})$ instead of $L_{k}(f)$.

Now the definitions of $m$-fold Riemann integrable functions $f\in RI(\mathbb{R}^{m})$
and their $m$-fold Riemann integrals $\int_{\mathbb{R}^{m}}f$ are
exactly the same as for $m=1$, and any properties of these functions
or their integrals that need to be proved can be proved exactly analogously.
In particular, in order to prove (\ref{eq:rwq}) (and similarly (\ref{eq:eeq}))
in this context, we note that each $I\in D_{k}(\mathbb{R}^{m})$ is
the disjoint union of $2^{m}$ dyadic $m$-intervals in $D_{k+1}(\mathbb{R}^{m})$,
which we may denote by $J_{1}(I)$, $J_{2}(I)$,...., $J_{2^{m}}(I)$
and that, analogously to (\ref{eq:blrp}), 
\[
2^{-(k+1)m}\sum_{n=1}^{2^{m}}\sup_{p\in J_{n}(I)}f(p)\le2^{-k\left(m+1\right)}\sum_{n=1}^{2^{m}}\sup_{p\in I}f(p)=2^{-km}\sup_{p\in I}f(p).
\]
We then sum over all $I\in D_{k}(\mathbb{R}^{m})$ to obtain (\ref{eq:rwq}). 
\begin{rem}
It should be mentioned that for $m=2$ and $f\in RI(\mathbb{R}^{2})$
the integral $\int_{\mathbb{R}^{2}}f$, is often called a \textbf{\textit{double
integral}}, and similarly, for $m=3$, and $f\in RI(\mathbb{R}^{3})$
the integral $\int_{\mathbb{R}^{3}}f$ is often called a \textbf{\textit{triple
integral}} or sometimes a \textbf{\textit{volume integral}}. The more
``traditional'' notation for these integrals is, respectively, $\iint_{E}f(x,y)dxdy$
or $\iint_{E}f(x,y)dA$ and $\iiint_{E}f(x,y,z)dxdydz$ or $\iiint_{E}f(x,y,z)dV$,
where $E$ is some subset of $\mathbb{R}^{2}$ or $\mathbb{R}^{3}$
respectively, such that $f$ vanishes on the complement of $E$ or
is simply not defined on that complement. For arbitrary $m$ the notation
for $\int_{\mathbb{R}^{m}}f$ can be, as for example in \cite{S},
something like $\int_{E}f(x^{1},x^{2},..,x^{m})dx^{1}dx^{2}...dx^{m}$.

Integrals of course have many applications and connections with various
problems in geometry, physics, etc. etc. Probably the main motivating
connections are that, for non-negative functions $f$, for $m=1$,
$2$ and $3$, the integral $\int_{\mathbb{R}^{m}}f$ gives us, respectively,
the area under the graph of $y=f(x)$, the volume under the graph
of $z=f(x,y)$ and the mass of a three dimensional region $E$ of
possibly varying density, where that density is given by $f(x,y,z)$
at each point $(x,y,z)\in E$.
\end{rem}

\begin{rem}
Dyadic intervals, squares, cubes and $m$-cubes have many other interesting
and important applications in mathematics apart from their usefulness
in developing integration theory. For some purposes, many authors
prefer to define them to be (products of) closed intervals instead
of semiclosed intervals.
\end{rem}
{\small }%
{\small \par}

{\small }%
{\small \par}
\begin{thm}
\label{thm:Monotone}Suppose that $f$ and $g$ are in $BB\left(\mathbb{R}^{m}\right)$
and satisfy $f\le g$. Then $U_{k}(f)\le U_{k}(g)$ and $L_{k}(f)\le L_{k}(g)$
for each $k\in\mathbb{N}$. Furthermore, if $f$ and $g$ are also
in $RI\left(\mathbb{R}^{m}\right)$, then $\int_{\mathbb{R}^{m}}f\le\int_{\mathbb{R}^{m}}g$.
\end{thm}

\noindent \textit{Proof.} This follows immediately from (\ref{eq:sp-mon}).
$\qed$
\begin{thm}
\label{thm:cubeok}For each $m\in\mathbb{N}$, each $k\in\mathbb{N}$
and each $I\in D_{k}(\mathbb{R}^{m})$ the function $\chi_{I}$ is
in $RI\left(\mathbb{R}^{m}\right)$ and $\int_{\mathbb{R}^{m}}\chi_{I}=2^{-km}$.
\end{thm}
\textit{Proof.} Obviously $\inf_{p\in I}\chi_{I}(p)=\sup_{p\in I}\chi_{I}(p)=1$
and, for every other $m$-dimensional ``dyadic cube'' $J$ in $D_{k}\left(\mathbb{R}^{m}\right)$
we have $\inf_{p\in J}\chi_{I}(p)=\sup_{p\in J}\chi_{I}(p)=0$ since
such $J$ is disjoint from $I$. Thus we obtain that $L_{k}(\chi_{I})=U_{k}(\chi_{I})=2^{-mk}$.
But, by general properties of the sequences $\left\{ L_{k}(f)\right\} _{k\in\mathbb{N}}$
and $\left\{ U_{k}\left(f\right)\right\} _{k\in\mathbb{N}}$ we know
that, for each $k'\ge k$ we have 
\begin{equation}
L_{k}(\chi_{I})\le L_{k'}(\chi_{I})\le U_{k'}\left(\chi_{I}\right)\le U_{k}(\chi_{I}).\label{eq:vvvglfd}
\end{equation}
It follows that all the inequalities in (\ref{eq:vvvglfd}) must be
equalities for each $k'\ge k$ and so we can pass to the limit and
obtain the desired conclusion. $\qed$ 
\begin{thm}
\label{thm:rlinearity}For all functions $f$ and $g$ in $RI\left(\mathbb{R}^{m}\right)$
and for all real constants $\alpha$ and $\beta$, the function $\alpha f+\beta g$
is an element of $RI\left(\mathbb{R}^{m}\right)$ and 
\[
\int_{\mathbb{R}^{m}}\left(\alpha f+\beta g\right)=\alpha\int_{\mathbb{R}^{m}}f+\beta\int_{\mathbb{R}^{m}}g.
\]

\end{thm}

\noindent \textit{Proof.} We proceed via several easy steps which
will be left as exercises in this preliminary version of this document. 

(i) The case where $\alpha=\beta=1$. (Do this by using (\ref{eq:sp-add})
and (\ref{eq:sp-add-inf}) to obtain that $U_{k}(f+g)\le U_{k}(f)+U_{k}(g)$
and $L_{k}(f)+L_{k}(g)\le L_{k}(f+g)$.)

(ii) The case where $\alpha\ge0$ and $\beta=0$. (Use (\ref{eq:sp01}).)

{} (iii) The case where $\alpha<0$ and $\beta=0$. (Use (\ref{eq:sp02}).)

{} (iv) The general case: By cases (ii) and %
{} (iii) we obtain that $\alpha f$ and $\beta g$ are both in $RI\left(\mathbb{R}^{m}\right)$
and that $\int_{\mathbb{R}^{m}}\alpha f=\alpha\int_{\mathbb{R}^{m}}f$
and $\int_{\mathbb{R}^{m}}\beta g=\beta\int_{\mathbb{R}^{m}}g$. Now
case (i), applied where the functions $f$ and $g$ are replaced by
$\alpha f$ and $\beta g$, tells us that $\alpha f+\beta g\in RI\left(\mathbb{R}^{m}\right)$
and, together with the previous cases, that 
\[
\int_{\mathbb{R}^{m}}\left(\alpha f+\beta g\right)=\int_{\mathbb{R}^{m}}\alpha f+\int_{\mathbb{R}^{m}}\beta g=\alpha\int_{\mathbb{R}^{m}}f+\beta\int_{\mathbb{R}^{m}}g.
\]
$\qed$

Here is one immediate consequence of Theorem \ref{thm:rlinearity}.
\begin{cor}
\label{cor:FinAditiv}Suppose that the sets $E_{1}$, $E_{2}$ ,....,
$E_{N}$ are pairwise disjoint subsets of $\mathbb{R}^{m}$ and that
$\chi_{E_{j}}\in RI(\mathbb{R}^{m})$ for $j=1,2...,N$. Then the
set $E:=\bigcup_{j=1}^{N}E_{j}$ satisfies $\chi_{E}\in RI(\mathbb{R}^{m})$
and 
\begin{equation}
\int_{\mathbb{R}^{m}}\chi_{E}=\sum_{j=1}^{N}\int_{\mathbb{R}^{m}}\chi_{E_{j}}.\label{eq:FinAdditiv}
\end{equation}

\end{cor}
\noindent \textit{Proof.} We simply use the fact that $\chi_{E}=\sum_{j=1}^{n}\chi_{E_{j}}$.
$\qed$
\begin{thm}
\label{thm:DyadicMultiCube} (i) Let $[a,b)$ be an interval such
that $a=n2^{-k}$ and $b=n'2^{-k}$ for some $n$ and $\, n'\in\mathbb{Z}$
and $k\in\mathbb{N}\cup\left\{ 0\right\} $. Then $\chi_{[a,b)}\in RI(\mathbb{R})$
and $\int_{\mathbb{R}}\chi_{[a,b)}=b-a$. 

(ii) More generally, suppose that $E=\prod_{j=1}^{m}[a_{j},b_{j})$
is an $m$-rectangle for some $m\in\mathbb{N}$, and that, for some
$k\in\mathbb{N\cup}\left\{ 0\right\} $, the endpoints of each of
the ``factor'' intervals $[a_{j},b_{j})$ are of the form $a_{j}=n_{j}2^{-k}$
and $b_{j}=n_{j}^{\prime}2^{-k}$ for some $n_{j}$ and $n_{_{j}}^{\prime}\in\mathbb{Z}$.
Then $\chi_{E}\in RI(\mathbb{R}^{m})$ and 
\begin{equation}
\int_{\mathbb{R}^{m}}\chi_{E}=\left|E\right|_{m}.\label{eq:VolOfMultiCube}
\end{equation}

\end{thm}
\noindent \textit{Proof.} For part (i) we observe that the interval
$[a,b)$ is the disjoint union $\bigcup$$_{\nu=n}^{n'-1}[\nu2^{-k},\left(\nu+1\right)2^{-k})$
and therefore $\chi_{[a,b)}=\sum_{\nu=n_{1}}^{n'-1}\chi_{[\nu2^{-k},\left(\nu+1\right)2^{-k})}$.
By Theorem \ref{thm:cubeok}, $\int_{\mathbb{R}}\chi_{[\nu2^{-k},\left(\nu+1\right)2^{-k})}=2^{-k}$
for each $\nu$, and so, by Theorem \ref{thm:rlinearity}, $\chi_{[a,b)}\in RI(\mathbb{R})$
and $\int_{\mathbb{R}}\chi_{[a,b)}$ is the sum of $n'-n$ terms each
equal to $2^{-k}$, i.e., it equals $\left(n'-n\right)2^{-k}=b-a$. 

For part (ii), for each $j\in\left\{ 1,2,...,m\right\} $ let $\mathcal{A}_{j}$
denote the collection of the pairwise disjoint intervals in $D_{k}(\mathbb{R})$
whose union is $[a_{j},b_{j})$. Then $E$ is the disjoint union of
all $m$-cubes $J$ of the form $J=\prod_{j=1}^{m}I_{j}$ where $I_{j}\in\mathcal{A}_{j}$
for each $j\in\left\{ 1,2,...,m\right\} $. Therefore $\chi_{E}$
is the sum of the characteristic functions of these cubes. There are
$\prod_{j=1}^{m}\left(n_{j}^{\prime}-n_{j}\right)$ such cubes, since
$\mathcal{A}_{j}$ consists of $n_{j}^{\prime}-n_{j}$ intervals for
each $j$. Thus, again using Theorem \ref{thm:cubeok} and then Theorem
\ref{thm:rlinearity}, we obtain that $\chi_{E}\in RI\left(\mathbb{R}^{m}\right)$
and $\int_{\mathbb{R}^{m}}\chi_{E}$ is the sum of $\prod_{j=1}^{m}\left(n_{j}^{\prime}-n_{j}\right)$
integrals of characteristic functions of cubes in $D_{k}\left(\mathbb{R}^{m}\right)$.
It follows that 
\[
\int_{\mathbb{R}^{m}}\chi_{E}=2^{-km}\prod_{j=1}^{m}\left(n_{j}^{\prime}-n_{j}\right)=\prod_{j=1}^{m}\left[\left(n_{j}^{\prime}-n_{j}\right)2^{-k}\right]=\prod_{j=1}^{m}(b_{j}-a_{j})=\left|E\right|_{m}
\]
as required. $\qed$ 

The following result is an immediate and obvious consequence of (\ref{eq:rwq}),
(\ref{eq:eeq}), (\ref{eq:burz}) and the definition of integrability. 
\begin{fact}
\label{fact:ec4ri}For each $f\in BB\left(\mathbb{R}^{m}\right)$
the sequence $U_{k}(f)-L_{k}(f)$ is non increasing and non-negative.
The function $f$ is in $RI\left(\mathbb{R}^{m}\right)$ if and only
if 
\begin{equation}
\lim_{k\to\infty}U_{k}(f)-L_{k}(f)=0\label{eq:isw}
\end{equation}
or, equivalently, 
\begin{equation}
\lim_{k\to\infty}\left(2^{-km}\sum_{I\in D_{k}(\mathbb{R}^{m})}\sup_{p,q\in I}\left|f(p)-f(q)\right|\right)=0.\label{eq:iswz}
\end{equation}

\end{fact}
\bigskip{}

What kind of functions in $BB(\mathbb{R}^{m})$ are integrable? We
will see that continuity is a sufficient condition for such functions
to be integrable. But we will also see that it is not necessary to
require such functions to be continuous absolutely everywhere. That
is to say, if the set of points where they are not continuous is ``very
small'' in some suitable sense, then such functions will still be
integrable. We will need several steps and auxiliary results to explain
and prove this. For now we will only present the first of those steps: 
\begin{thm}
\label{thm:NeedUnifCty}Suppose that $I\in D_{k_{0}}(\mathbb{R}^{m})$
for some $k_{0}\in\mathbb{N}\cup\left\{ 0\right\} $ and some $m\in\mathbb{N}$.
Suppose that the function $f:\mathbb{R}^{m}\to\mathbb{R}$ satisfies
$f(p)=0$ for all $p\in\mathbb{R}^{m}\setminus I$ and is uniformly
continuous on $I$. Then $f\in RI\left(\mathbb{R}^{m}\right)$.
\end{thm}
\noindent \textit{Proof.} Since the set $I$ is contained in a compact
set, and $f$ is uniformly continuous on $I$ it follows that $f$
is bounded on $I$ (because $f$ maps Cauchy sequences in $I$ to
Cauchy sequences in $\mathbb{R}$). The uniform continuity of $f$
on $I$ implies that for each $\epsilon>0$ there exists a number
$\delta(\epsilon)>0$ such that, whenever two points $p=\left(p_{1},p_{2},...,p_{m}\right)$
and $q=\left(q_{1},q_{2},....,q_{m}\right)$ in $I$ satisfy $\max\left\{ \left|p_{j}-q_{j}\right|:j=1,2,...,m\right\} <\delta(\epsilon)$,
then $\left|f(p)-f(q)\right|<\epsilon$. 

So, given any $\epsilon>0$, we choose an integer $k\ge k_{0}$ such
that $2^{-k}<\delta(\epsilon)$. For each $m$-cube $J\in D_{k}(\mathbb{R}^{m})$
which is contained in $I$ and for each pair of points $p,q\in J$
we have $\max\left\{ \left|p_{j}-q_{j}\right|:j=1,2,...,m\right\} <2^{-k}<\delta(\epsilon)$
and therefore $\left|f(p)-f(q)\right|<\epsilon$. This means that
$\sup_{p,q\in J}\left|f(p)-f(q)\right|\le\epsilon$. Let us now see
that this in turn implies that 
\begin{equation}
2^{-km}\sum_{J\in D_{k}(\mathbb{R}^{m})}\sup_{p,q\in J}\left|f(p)-f(q)\right|\le2^{-k_{0}m}\epsilon.\label{eq:wwst}
\end{equation}
Since $f$ vanishes on $\mathbb{R}^{m}\setminus I$ the sum on the
left side of (\ref{eq:wwst}) reduces to the sum over all $J\in D_{k}(\mathbb{R}^{m})$
which are contained in $I$. For each $n\ge0$, each $m$-cube in
$D_{n}(\mathbb{R}^{m})$ is the disjoint union of $2^{m}$ $m$-cubes
in $D_{n+1}\left(\mathbb{R}^{m}\right)$. By repeated application
of this fact we see that the $m$-cube $I$ is the disjoint union
of $\left(2^{m}\right)^{k-k_{0}}$ $m$-cubes in $D_{k}\left(\mathbb{R}^{m}\right)$.
Thus there are $2^{m(k-k_{0})}$ non-zero terms in the above-mentioned
sum, each dominated by $\epsilon$ and the left side of (\ref{eq:wwst})
therefore does not exceed $2^{-km}\times2^{m(k-k_{0})}\epsilon=2^{-k_{0}m}\epsilon$.
Since $\epsilon$ is arbitrary and $2^{-k_{0}m}\le1$ we have established
(\ref{eq:iswz}) which shows that $f\in RI\left(\mathbb{R}^{2}\right)$.
$\qed$
\begin{rem}
It is easy to use Theorem \ref{thm:NeedUnifCty} to show that any
function $f\in(\mathbb{R}^{m})$ which is continuous (and therefore
uniformly continuous) is also in $RI\left(\mathbb{R}^{m}\right)$.
But we do not need that result at this stage and in any case we will
later prove a stronger result (namely Theorem \ref{thm:CtsAE}). \medskip{}

\end{rem}
The next two theorems give us several ways of combining functions
which are known to be integrable in order to obtain new integrable
functions.
\begin{thm}
\label{thm:ComposedFn} Let $\phi:\mathbb{R}\to\mathbb{R}$ be a function
whose restriction to each bounded interval is a Lipschitz function,
i.e., for each $r>0$ there exists a constant $\lambda_{r}$ such
that $\left|\phi(x)-\phi(x')\right|\le\lambda_{r}\left|x-x'\right|$
for each pair of points $x$ and $x'$ in $\left[-r,r\right]$. Then
the composed function $\phi\circ f$ is in $RI\left(\mathbb{R}^{n}\right)$
whenever $f\in RI\left(\mathbb{R}^{n}\right)$.
\end{thm}
\noindent \textit{Proof.} Given $f\in RI\left(\mathbb{R}^{n}\right)$
we let $r=\sup_{p\in\mathbb{R}^{n}}\left|f(p)\right|$. Then, for
each $k\in\mathbb{N}$ and for each $I\in D_{k}\left(\mathbb{R}^{n}\right)$,
and for each $p$ and $q$ in $I$, we have $f(p),f(q)\in\left[-r,r\right]$
and therefore $\left|\phi\circ f(p)-\phi\circ f(q)\right|\le\lambda_{r}\left|f(p)-f(q))\right|.$
It follows that
\begin{equation}
2^{-km}\sum_{I\in D_{k}(\mathbb{R}^{m})}\sup_{p,q\in I}\left|\phi\circ f(p)-\phi\circ f(q)\right|\le\lambda_{r}2^{-km}\sum_{I\in D_{k}(\mathbb{R}^{m})}\sup_{p,q\in I}\left|f(p)-f(q)\right|.\label{eq:korq}
\end{equation}
Since, by Fact \ref{fact:ec4ri}, the right side of (\ref{eq:korq})
tends to $0$ as $k\to\infty$, the same is true for the left side,
and consequently, again by Fact \ref{fact:ec4ri}, $\phi\circ f\in RI(\mathbb{R}^{n})$.
$\qed$
\begin{rem}
A more elaborate proof would show that the result of Theorem \ref{thm:ComposedFn}
also holds when $\phi$ is merely continuous.\end{rem}
\begin{thm}
\label{thm:CorollOfComposedFn}For each $f$ and $g$ in $RI(\mathbb{R}^{n})$
the functions $fg$, $\left|f\right|$, $f_{+}:=\max\left\{ f,0\right\} $,
$f_{-}:=\max\left\{ -f,0\right\} $, $\max\left\{ f,g\right\} $ and
$\min\left\{ f,g\right\} $ are also in $RI\left(\mathbb{R}^{n}\right)$.
\end{thm}
\noindent \textit{Proof.} Since $\phi(t)=t^{2}$ satisfies the hypotheses
of Thorem \ref{thm:ComposedFn} we obtain that $h^{2}\in RI(\mathbb{R}^{n})$
for all $h\in RI(\mathbb{R}^{n})$. Since $fg=\frac{1}{4}\left((f+g)^{2}-\left(f-g\right)^{2}\right)$
we deduce that $fg\in RI\left(\mathbb{R}^{n}\right)$ for all $f,g\in RI\left(\mathbb{R}^{n}\right)$. 

Since the function $\phi(t)=\left|t\right|$ satisfies the hypotheses
of Theorem \ref{thm:ComposedFn} we obtain that $\left|f\right|=\phi\circ f\in RI\left(\mathbb{R}^{n}\right)$.
Since the functions $f_{+}$ and $f_{-}$ are given by the formulae
$f_{+}=\left(\left|f\right|+f\right)/2$ and $f_{-}=\left(\left|f\right|-f\right)/2$
they too are in $RI\left(\mathbb{R}^{n}\right)$. Finally, we reach
the desired conclusions for the remaining functions, by observing
that $\max\left\{ f,g\right\} =\max\left\{ f-g,0\right\} +g$ and
$\min\left\{ f,g\right\} =-\max\left\{ -f,-g\right\} $. $\qed$

\medskip{}

Next we present another condition for integrability which will be
a convenient tool in the subsequent development of the theory. 
\begin{thm}
\label{thm:sandwich} Suppose that $f\in BB(\mathbb{R}^{m})$ and
that for each $\epsilon>0$ there exist a pair of functions $g_{\epsilon},h_{\epsilon}\in RI(\mathbb{R}^{m})$
such that 
\begin{equation}
g_{\epsilon}\le f\le h_{\epsilon}\label{eq:btiv}
\end{equation}
 and 
\begin{equation}
\int_{\mathbb{R}^{m}}\left(h_{\epsilon}-g_{\epsilon}\right)\le\epsilon.\label{eq:bciw}
\end{equation}
Then $f\in RI(\mathbb{R}^{m})$ and 
\begin{equation}
\int_{\mathbb{R}^{m}}f=\lim_{\epsilon\searrow0}\int_{\mathbb{R}^{m}}g_{\epsilon}=\lim_{\epsilon\searrow0}\int_{\mathbb{R}^{m}}h_{\epsilon}.\label{eq:FeGepHep}
\end{equation}

Consequently, if the number $\alpha$ satisfies $\int_{\mathbb{R}^{m}}g_{\epsilon}\le\alpha\le\int_{\mathbb{R}^{m}}h_{\epsilon}$
for each $\epsilon>0$ and the functions $g_{\epsilon}$ and $h_{\epsilon}$
with the above properties, then $\alpha=\int_{\mathbb{R}^{m}}f$.
\end{thm}
\textit{Proof.} Given any fixed $\epsilon>0$ there exist functions
$h_{\epsilon/2}$ and $g_{\epsilon/2}$ depending on $\epsilon$ such
that, for each $k\in\mathbb{N}$, we can use (\ref{eq:btiv}), Theorem
\ref{thm:Monotone} and then (\ref{eq:bciw}) to obtain the following
inequalities.

\begin{eqnarray*}
 &  & U_{k}(f)-L_{k}(f)\\
 & \le & U_{k}(h_{\epsilon/2})-L_{k}(g_{\epsilon/2})\\
 & = & \left[U_{k}(h_{\epsilon/2})-\int_{\mathbb{R}^{m}}h_{\epsilon/2}\right]+\left[\int_{\mathbb{R}^{m}}g_{\epsilon/2}-L_{k}(g_{\epsilon/2})\right]+\int_{\mathbb{R}^{m}}\left(h_{\epsilon/2}-g_{\epsilon/2}\right)\\
 & \le & \left[U_{k}(h_{\epsilon/2})-\int_{\mathbb{R}^{m}}h_{\epsilon/2}\right]+\left[\int_{\mathbb{R}^{m}}g_{\epsilon/2}-L_{k}(g_{\epsilon/2})\right]+\frac{\epsilon}{2}.
\end{eqnarray*}

Since $h_{\epsilon/2}$ and $g_{\epsilon/2}$ are both in $RI\left(\mathbb{R}^{m}\right)$
there exists an integer $k\left(\epsilon\right)$ such that, for all
$k\ge k(\epsilon)$ each of the expressions in the square brackets
above does not exceed $\epsilon/4$. Thus we have shown that for each
$k\ge k(\epsilon)$ we also have 
\[
U_{k}(f)-L_{k}(f)\le\frac{\epsilon}{2}+\frac{\epsilon}{4}+\frac{\epsilon}{4}=\epsilon.
\]
 This shows that $f\in RI(\mathbb{R}^{m})$. We can now assert (by
Theorem\ref{thm:Monotone}) that $\int_{\mathbb{R}^{m}}g_{\epsilon}\le\int_{\mathbb{R}^{m}}f\le\int_{\mathbb{R}^{m}}h_{\epsilon}$
for each $\epsilon>0$. This implies that $\left|\int_{\mathbb{R}^{m}}f-\int_{\mathbb{R}^{m}}g_{\epsilon}\right|=\int_{\mathbb{R}^{m}}f-\int_{\mathbb{R}^{m}}g_{\epsilon}\le\int_{\mathbb{R}^{m}}h_{\epsilon}-\int_{\mathbb{R}^{m}}g_{\epsilon}<\epsilon$.
Similarly $\left|\int_{\mathbb{R}^{m}}f-\int_{\mathbb{R}^{m}}h_{\epsilon}\right|=\int_{\mathbb{R}^{m}}h_{\epsilon}-\int_{\mathbb{R}^{m}}f\le\int_{\mathbb{R}^{m}}h_{\epsilon}-\int_{\mathbb{R}^{m}}g_{\epsilon}<\epsilon$.
These inequalities imply (\ref{eq:FeGepHep}) and so complete the
proof. $\qed$.

The next definition and proposition provide us with a convenient variant
of the preceding result.
\begin{defn}
A function $f\in BB(\mathbb{R}^{n})$ will be said to be a \textbf{\textit{dyadic
step function}} or, more explicitly, a \textbf{\textit{dyadic step
function of order $k$}}, if for some $k\in\mathbb{N}$ it can be
written in the form $f=\sum_{I\in D_{k}(\mathbb{R}^{n})}\alpha_{I}\chi_{I}$
for some real numbers $\alpha_{I}$. The fact that $f\in BB(\mathbb{R}^{n})$
ensures that only finitely many of the numbers $\alpha_{I}$ can be
non-zero. \end{defn}
\begin{prop}
\label{prop:DyadicSandwich}Let $f\in BB(\mathbb{R}^{m})$. Then $f\in RI(\mathbb{R}^{m})$
if and only if there exist two sequences $\left\{ g_{k}\right\} _{k\in\mathbb{N}}$
and $\left\{ h_{k}\right\} _{k\in\mathbb{N}}$ of dyadic step functions
such that 
\begin{equation}
g_{k}\le f\le h_{k}\mbox{ for each }k\in\mathbb{N}\label{eq:yhstzm}
\end{equation}
and 
\begin{equation}
\lim_{k\to\infty}\int_{\mathbb{R}^{m}}\left(h_{k}-g_{k}\right)=0.\label{eq:hbisl}
\end{equation}

When $f\in RI(\mathbb{R}^{m})$ any sequences $\left\{ g_{k}\right\} _{k\in\mathbb{N}}$
and $\left\{ h_{k}\right\} _{k\in\mathbb{N}}$ with the above properties
satisfy 
\[
\lim_{k\to\infty}\int_{\mathbb{R}^{m}}g_{k}=\lim_{k\to\infty}\int_{\mathbb{R}^{m}}h_{k}=\int_{\mathbb{R}^{m}}f
\]
 and any number $\alpha$ satisfying $\int_{\mathbb{R}^{m}}g_{k}\le\alpha\le\int_{\mathbb{R}^{m}}h_{k}$
for all $k\in\mathbb{N}$ must equal $\int_{\mathbb{R}^{m}}f$.
\end{prop}

\noindent \textit{Proof.} Suppose first that $f\in RI(\mathbb{R}^{n})$.
For each $k\in\mathbb{N}$ let $g_{k}=\sum_{I\in D_{k}(\mathbb{R}^{n})}\inf_{p\in I_{k}}f(p)\chi_{I}$
and $h_{k}=\sum_{I\in D_{k}(\mathbb{R}^{n})}\sup_{p\in I_{k}}f(p)\chi_{I}$.
These functions clearly satisfy (\ref{eq:yhstzm}). Since both of
these sums have only finitely many non-zero terms we can use Theorem
\ref{thm:rlinearity} and then Theorem \ref{thm:cubeok} to obtain
that 
\[
\int_{\mathbb{R}^{n}}g_{k}=\sum_{I\in D_{k}(\mathbb{R}^{n})}\inf_{p\in I_{k}}f(p)\int_{\mathbb{R}^{n}}\chi_{I_{k}}=\sum_{I\in D_{k}(\mathbb{R}^{n})}\inf_{p\in I_{k}}f(p)\int_{\mathbb{R}^{n}}\chi_{I}=2^{-kn}\sum_{I\in D_{k}(\mathbb{R}^{n})}\inf_{p\in I_{k}}f(p)=L_{k}(f).
\]
Similarly $\int_{\mathbb{R}^{n}}h_{k}=U_{k}(f)$ and we see that (\ref{eq:hbisl})
is simply (\ref{eq:isw}), which holds by Fact \ref{fact:ec4ri}. 

For the reverse implication we simply apply Theorem \ref{thm:sandwich}. 

The claims in the final sentence of the theorem are immediate consequences
of (\ref{eq:yhstzm}) and (\ref{eq:hbisl}). $\qed$

The following two results will enable us to extend the facts that
$\chi_{I}\in RI(\mathbb{R}^{m})$ and $\int_{\mathbb{R}^{m}}\chi_{I}=\left|I\right|_{m}$
which we obtained very easily (see Theorem \ref{thm:cubeok}) in the
case where $I$ is a dyadic interval, dyadic square, dyadic cube or
dyadic $m$-cube, to the case where $I$ is any $m$-rectangle. 
\begin{prop}
\label{prop:ArbitraryRectangle}Let $E$ be an arbitrary $m$-rectangle
with non-empty interior. Then there exist two sequences $\left\{ A_{n}\right\} _{n\in\mathbb{N}}$
and $\left\{ B_{n}\right\} _{n\in\mathbb{N}}$ of $m$-rectangles
with the following properties:
\begin{equation}
(i)\qquad\qquad\,\,\,\,\,\,\,\,\lim_{n\to\infty}\left|A\right|_{m}=\lim_{n\to\infty}\left|B_{n}\right|_{m}=\left|E\right|_{m}.\,\,\,\,\,\,\,\,\qquad\qquad\label{eq:AeBeE}
\end{equation}
(ii) For each sufficiently large $n$, $A_{n}$ and $B_{n}$ satisfy
\begin{equation}
A_{n}\subset\overline{A_{n}}\subset E^{\circ}\subset E\subset\overline{E}\subset B_{n}^{\circ}\subset B_{n},\label{eq:Inclusions}
\end{equation}
(iii) The $m$-rectangles $A_{n}$ and $B_{n}$ are both finite unions
of cubes in $D_{n}\left(\mathbb{R}^{m}\right)$, the functions $\chi_{A_{n}}$
and $\chi_{B_{n}}$ are both in $RI\left(\mathbb{R}^{m}\right)$ and
their integrals satisfy 
\[
\int_{\mathbb{R}^{m}}\chi_{A_{n}}=\left|A_{n}\right|_{m}\mbox{ and }\int_{\mathbb{R}^{m}}\chi_{B_{n}}=\left|B_{n}\right|_{m}.
\]

\end{prop}
\noindent \textit{Proof.} For each $t\in\mathbb{R}$ let $I_{n.t}$
be the unique interval in $D_{n}(\mathbb{R})$ which contains $t$.
Let $\gamma_{n}(t)$ be the left endpoint of the interval in $D_{n}(\mathbb{R})$
which is adjacent to $I_{n,t}$ on its left side, and let $\delta_{n}(t)$
be the right endpoint of $I_{n,t}$. Then 
\begin{equation}
\gamma_{n}(t)<t<\delta_{n}(t)\label{eq:GamLtDel}
\end{equation}
 and $\delta_{n}(t)-\gamma_{n}(t)=2^{1-n}$. We observe that, for
any $a$ and $b$ such that $a<b$, we have 
\begin{equation}
\gamma_{n}(b)-\delta_{n}(a)=\delta_{n}(b)-2^{1-n}-(\gamma_{n}(a)+2^{1-n})=\delta_{n}(b)-\gamma_{n}(a)-2^{2-n}>b-a-2^{2-n}.\label{eq:GamDelStuff}
\end{equation}
 The open $m$-rectangle $E^{\circ}$ is of course of the form $\prod_{j=1}^{m}(a_{j},b_{j})$.
For all sufficiently large $n$ we set $A_{n}=\prod_{j=1}^{m}[\delta_{n}(a_{j}),\gamma_{n}(b_{j}))$
and $B_{n}=\prod_{j=1}^{m}[\gamma_{n}(a_{j}),\delta_{n}(b_{j}))$.
``Sufficiently large'' means that we need to have 
\begin{equation}
\delta_{n}(a_{j})<\gamma_{n}(b_{j})\mbox{ for all }j\in\{1,2,...,m\},\label{eq:szum}
\end{equation}
and (\ref{eq:GamDelStuff}) shows that this can be guaranteed by choosing
$n$ large enough to satisfy 
\[
2^{2-n}<\min\left\{ (b_{j}-a_{j}):j=1,2,...,m\right\} .
\]
The inclusions $\overline{A}_{n}\subset E^{\circ}$ and $\overline{E}\subset B_{n}^{\circ}$
follow immediately from (\ref{eq:GamLtDel}) and (\ref{eq:szum}).
Since $\lim_{n\to\infty}\gamma_{n}(t)=t$ and $\lim_{n\to\infty}\delta_{n}(t)=t$
for all $t\in\mathbb{R}$, we see that 
\[
\lim_{n\to\infty}\left|A\right|_{m}=\lim_{n\to\infty}\left|B_{n}\right|_{m}=\prod_{j=1}^{m}(b_{j}-a_{j})=\left|E\right|_{m}
\]
establishing (\ref{eq:AeBeE}). Finally, the conclusions of part (iii)
of the proposition follow immediately from Theorem \ref{thm:DyadicMultiCube}.
$\qed$
\begin{cor}
\label{cor:AllRectanglesMeasurable}Let $E$ be an arbitrary $m$-rectangle.
Then the functions $\chi_{E}$, $\chi_{E^{\circ}}$ and $\chi_{\overline{E}}$
are all in $RI\left(\mathbb{R}^{m}\right)$ and 
\begin{equation}
\int_{\mathbb{R}^{m}}\chi_{E}=\int_{\mathbb{R}^{m}}\chi_{E^{\circ}}=\int_{\mathbb{R}^{m}}\chi_{\overline{E}}=\left|E\right|_{m}.\label{eq:AllEqual}
\end{equation}

\end{cor}
\noindent \textit{Proof.} Let us first consider the case where $E$
has non-empty interior. Then, for each $\epsilon>0$ we can obtain
functions $g_{\epsilon}$ and $h_{\epsilon}$ in $RI\left(\mathbb{R}^{m}\right)$
which satisfy $g_{\epsilon}\le\chi_{E^{\circ}}\le\chi_{E}\le\chi_{\overline{E}}\le h_{\epsilon}$
and $\int_{\mathbb{R}^{m}}(h_{\epsilon}-g_{\epsilon})<\epsilon$ by
choosing $g_{\epsilon}=\chi_{A_{n}}$ and $h_{\epsilon}=\chi_{B_{n}}$
where $A_{n}$ and $B_{n}$ are elements of the sequence constructed
in Proposition \ref{prop:ArbitraryRectangle}, with $n$ chosen sufficiently
large to ensure that $\left|B_{n}\right|_{m}-\left|A_{n}\right|_{m}<\epsilon$.
The required result follows by applying Theorem \ref{thm:sandwich}.

It remains to consider the case where $E^{\circ}$ is empty. Then
$\overline{E}=\prod_{j=1}^{m}[a_{j},b_{j}]$ where $a_{j}=b_{j}$
for at least one value of $j$. Given any $\epsilon>0$ we choose
$\rho>0$ sufficiently small so that $E(\rho):=\prod_{j=1}^{m}[a_{j},b_{j}+\rho]$
satisfies $\left|E(\rho)\right|_{m}<\epsilon/2$. Since $E(\rho)$
has non-empty interior we can apply Proposition \ref{prop:ArbitraryRectangle}
to $E(\rho)$ to obtain $A_{n}$ and $B_{n}$ of the form stated in
that proposition, such that $A_{n}\subset\overline{E(\rho)}=E(\rho)\subset B_{n}$
and with $n$ chosen sufficiently large to ensure that $\left|B_{n}\right|_{m}\le\left|A_{n}\right|_{m}+\epsilon/2<\left|E(\rho)\right|+\epsilon/2<\epsilon$.
Now choose $g_{\epsilon}=0$ and $h_{\epsilon}=\chi_{B_{n}}$. We
have $0=g_{\epsilon}=\chi_{E^{\circ}}\le\chi_{E}\le\chi_{\overline{E}}\le\chi_{E\left(\rho\right)}\le h_{\epsilon}$.
So Theorem \ref{thm:sandwich}, again shows that $\chi_{E^{\circ}}$,
$\chi_{E}$ and $\chi_{\overline{E}}$ are all in $RI(\mathbb{R}^{m})$
and that the integrals of all these functions are dominated by  $\int_{\mathbb{R}^{n}}h_{\epsilon}<\epsilon$
for all $\epsilon$. Therefore, also in this case, we obtain (\ref{eq:AllEqual})
where here of course $\left|E\right|_{m}=0$.$\qed$

We are now ready to present the following stronger variant of Theorem
\ref{thm:NeedUnifCty}.
\begin{thm}
\label{thm:CtySuffices}Suppose that $I\in D_{k_{0}}(\mathbb{R}^{m})$
for some $k_{0}\in\mathbb{N}\cup\left\{ 0\right\} $ and some $m\in\mathbb{N}$.
Suppose that the function $f:\mathbb{R}^{m}\to\mathbb{R}$ satisfies
$f(p)=0$ for all $p\in\mathbb{R}^{m}\setminus I$ and is bounded
and continuous on $I$. Then $f\in RI\left(\mathbb{R}^{m}\right)$.
\end{thm}
\noindent \textit{Proof.} Let $R=\sup_{p\in\mathbb{R}^{m}}\left|f(p)\right|=\sup_{p\in I}\left|f(p)\right|$.
For each $\epsilon>0$ we apply Proposition \ref{prop:ArbitraryRectangle}
to $I$ to obtain an $m$-rectangle $A_{n}$ which is finite collection
of pairwise disjoint dyadic cubes $J_{1}$, $J_{2}$,....., $J_{N}$
in $D_{n}\left(\mathbb{R}^{m}\right)$ which are all contained in
a closed subset of $I$, namely $\overline{A_{n}}$, and where (by
using (\ref{eq:AeBeE})) the integer $n\ge k_{0}$ has been chosen
sufficiently large to ensure that $\left|I\right|_{m}-\left|A_{n}\right|_{m}=\left|I\right|_{m}-\sum_{k=1}^{N}\left|J_{k}\right|_{m}<\epsilon/2R$.
Since $I$ is a cube in $D_{k_{0}}(\mathbb{R}^{m})$ it is the disjoint
union of $M=2^{m(n-k_{0})}$ cubes in $D_{n}\left(\mathbb{R}^{m}\right)$.
So we can write $I=\bigcup_{k=1}^{M}J_{k}$, where , for $k=1$ to
$N$ the cubes $J_{k}$ are from the above-mentioned collection, and
the additional cubes $J_{k}$ for $k=N+1$ to $M$ satisfy 
\[
\sum_{k=N+1}^{M}\left|J_{k}\right|_{m}=\left|I\right|_{m}-\sum_{k=1}^{N}\left|J_{k}\right|_{m}<\epsilon/2R.
\]
Now let $g_{\epsilon}=\sum_{k=1}^{N}f\chi_{J_{k}}-R\sum_{k=N+1}^{M}\chi_{J_{k}}$
and $h_{\epsilon}=\sum_{k=1}^{N}f\chi_{J_{k}}+R\sum_{k=N+1}^{M}\chi_{J_{k}}$.
For $k=1$ to $N$ the function $f$ is uniformly continuous on the
cube $J_{k}$, since $J_{k}$ is contained in a closed (and obviously
bounded) subset of $I$. So Theorem \ref{thm:NeedUnifCty} gives us
that $f\chi_{J_{k}}\in RI\left(\mathbb{R}^{m}\right)$. We deduce
that $g_{\epsilon}$ and $h_{\epsilon}$ are both finite sums of functions
in $RI(\mathbb{R}^{m})$ and therefore themselves in $RI(\mathbb{R}^{m})$.
Our choice of $R$ guarantees that $g_{\epsilon}\le f\le h_{\epsilon}$.
Furthermore, 
\[
0\le\int_{\mathbb{R}^{m}}(h_{\epsilon}-g_{\epsilon})=\int_{\mathbb{R}^{m}}2R\left(\sum_{n=N+1}^{M}\chi_{J_{n}}\right)=2R\sum_{n=N+1}^{M}\left|J_{n}\right|_{m}<\epsilon.
\]
Since $\epsilon$ can be chosen arbitrarily small, we can apply Theorem
\ref{thm:sandwich} to obtain that $f\in RI\left(\mathbb{R}^{m}\right)$.
$\qed$

The previous theorem brings us quite close to a quite general result
about integrability of functions which are continuous everywhere except
on a ``very small'' set. Before presenting that result we have to
discuss what we mean here by ``very small''. For the purposes of
this document we shall introduce the following informal terminology.
\begin{defn}
\label{def:VerySmall}Let $E$ be a subset of $\mathbb{R}^{m}$. We
shall say that $E$ is a \textbf{\textit{very small subset of}} $\mathbb{R}^{m}$
if, for each $\epsilon>0$ there exists a finite collection $J_{1}$,
$J_{2}$, ...., $J_{N}$ of dyadic $m$-cubes such that $E\subset\bigcup_{n=1}^{N}J_{n}$
and $\sum_{n=1}^{N}\left|J_{n}\right|_{m}<\epsilon$. 
\end{defn}

In more standard terminology very small subsets of $\mathbb{R}^{m}$
would be referred to as sets which have zero $m$-dimensional Jordan
measure. We will say something about Jordan measure later in this
document.

There are several equivalent ways of expressing the fact that a set
$E$ is a very small subset of $\mathbb{R}^{m}$. We shall now mention
a few of them:

Since each cube in $D_{k}(\mathbb{R}^{m})$ is a finite union of disjoint
cubes in $D_{k'}(\mathbb{R}^{m})$ for any $k'>k$, we see that we
can assume that all the cubes $J_{1},J_{2},....,J_{N}$ in any finite
collection appearing in Definition \ref{def:VerySmall} can be assumed
to all be in the same class $D_{k}(\mathbb{R}^{m})$ for some suitably
large $k$. Then these cubes can also of course be chosen to be all
different, i.e., pairwise disjoint. Thus the condition $\sum_{n=1}^{N}\left|J_{n}\right|_{m}<\epsilon$
implies that $U_{k}\left(\chi_{E}\right)<\epsilon$ for the above
choice of $k$, and is in fact equivalent to this last inequality.
From this we see that $E$ is a very small subset of $\mathbb{R}^{m}$
if and only if $\lim_{k\to\infty}U_{k}(\chi_{E})=0$. This in turn
is equivalent to the condition that $\chi_{E}\in RI(\mathbb{R}^{m})$
and satisfies $\int_{\mathbb{R}^{m}}\chi_{E}=0$. 

Using Theorem \ref{thm:sandwich} we can immediately see that another
equivalent condition for $E$ to be a very small subset of $\mathbb{R}^{m}$
is that, for each $\epsilon>0$, there exists a function $h_{\epsilon}$
which satisfies 
\begin{equation}
\chi_{E}\le h_{\epsilon}\mbox{ and }\int_{\mathbb{R}^{m}}h_{\epsilon}<\epsilon.\label{eq:cohe}
\end{equation}

Obviously every finite subset of $\mathbb{R}^{m}$ is a very small
subset of $\mathbb{R}^{m}$, but much larger subsets of $\mathbb{R}^{m}$
are also very small. For example, as we shall now see, the graph of
a continuous function of one variable on a closed bounded interval
is a very small subset of $\mathbb{R}^{2}$. Roughly speaking, this
means that this graph has zero ``area'', if we can decide what we
mean exactly by area.
\begin{thm}
\label{thm:GraphIsVerySmall}Let $\left[\alpha,\beta\right]$ be a
bounded closed interval and suppose that $\phi:\left[\alpha,\beta\right]\to\mathbb{R}$
is continuous. Then the set $E=\left\{ (x,y):\alpha\le x\le\beta,\, y=\phi(x)\right\} $
is a very small subset of $\mathbb{R}^{2}$. More generally, let $H$
be a bounded subset of $\mathbb{R}^{m}$ and let $\phi:H\to\mathbb{R}$
be a uniformly continuous function on $H$. Then the set 
\[
E=\left\{ (x_{1},x_{2},....,x_{m},x_{m+1})\in\mathbb{R}^{m+1}:(x_{1},x_{2},....,x_{m})\in H,\, x_{m+1}=\phi(x_{1},x_{2},....,x_{m})\right\} 
\]
is a very small subset of $\mathbb{R}^{m+1}$.
\end{thm}
\noindent \textit{Proof.} Please draw a picture which should show
you that, because of the uniform continuity of $\phi$ it is possible,
for each $\epsilon>0$, to contain $E$, the graph of $\phi$, in
a finite collection of sufficiently narrow rectangles each having
height less that $\epsilon/\left(\beta-\alpha\right)$. The sum of
the areas of these rectangles is less than $\epsilon$. So it is clear
intuitively that $E$ is ``very small'' in the required sense. 

Of course the first statement of the theorem is a special case of
the second statement. So let us just prove the second statement. Even
though the notation may be a little more complicated, the simple idea
of the proof will be the same as for the intuitive proof that you
were just encouraged to give for the first statement. 

Since $H$ is bounded, it can be contained in the set $H_{1}=\prod_{j=1}^{m}[-N,N)$
for some suitably large integer $N$. 

The uniform continuity of $\phi$ on $H$ implies that, for each $\epsilon>0$,
there exists some number $\delta(\epsilon)>0$ such that $\left|\phi(p)-\phi(q)\right|<(2N)^{-m}\epsilon$
for every pair of points $p$ and $q$ in $H$ which satisfy $\max\left\{ \left|p_{j}-q_{j}\right|:j=1,2,...,m\right\} <\delta(\epsilon)$.
Let us choose $k\in\mathbb{N}$ so that $2^{-k}<\delta(\epsilon)$.
The set $H_{1}$ is the disjoint union of all cubes in $D_{k}\left(\mathbb{R}^{m}\right)$
which it contains. Let $\mathcal{A}$ be the collection of all cubes
in $D_{k}(\mathbb{R}^{m})$ which have non-empty intersection with
$H$ and are therefore contained in $H_{1}$. We observe (using \ref{cor:AllRectanglesMeasurable}
and Corollary \ref{cor:FinAditiv}) that 
\begin{equation}
\sum_{I\in\mathcal{A}}\left|I\right|_{m}\le\sum_{I\in D_{k}(\mathbb{R}^{m}),I\subset H_{1}}\left|I\right|_{m}=\left|H_{1}\right|_{m}=(2N)^{m}.\label{eq:klomp}
\end{equation}

For each $I\in\mathcal{A}$, and for each $p$ and $q$ in $H\cap I$,
the facts that $I\in D_{k}(\mathbb{R}^{m})$ and $2^{-k}<\delta$
ensure that $\left|\phi(p)-\phi(q)\right|<(2N)^{-m}\epsilon$. So
(cf.~(\ref{eq:burz})), $\sup_{p\in H\cap J}\phi(p)-\inf_{p\in H\cap J}\phi(p)\le(2N)^{-m}\epsilon$.
Thus the $(m+1)$-rectangle 
\begin{eqnarray*}
J_{I} & := & I\times\left[\inf_{p\in H\cap J}\phi(p),\sup_{p\in H\cap J}\phi(p)\right]\\
 & = & \left\{ (x_{1},x_{2},....,x_{m},x_{m+1})\in\mathbb{R}^{m+1}:(x_{1},x_{2},....,x_{m})\in I,\,\inf_{p\in H\cap J}\phi(p)\le x_{m+1}\le\sup_{p\in H\cap J}\phi(p)\right\} 
\end{eqnarray*}
has $\left(m+1\right)$-volume $\left|J_{I}\right|_{m+1}<\left|I\right|_{m}(2N)^{-m}\epsilon$.
This rectangle contains all points $(x_{1},x_{2},....,x_{m},x_{m+1})\in E$
for which $(x_{1},x_{2},....,x_{m})\in H\cap I$ . Therefore $E$
is contained in $\bigcup_{I\in\mathcal{A}}J_{I}$ which implies that
\[
g_{\epsilon}:=0\le\chi_{E}\le h_{\epsilon}:=\sum_{I\in\mathcal{A}}\chi_{J_{I}}
\]
and 
\begin{eqnarray*}
\int_{_{\mathbb{R}^{m+1}}}(h_{\epsilon}-g_{\epsilon}) & = & \int_{_{\mathbb{R}^{m+1}}}h_{\epsilon}\le\sum_{I\in\mathcal{A}}\int_{\mathbb{R}^{m+1}}\chi_{J_{I}}=\sum_{I\in\mathcal{A}}\left|J_{I}\right|_{m+1}<(2N)^{-m}\epsilon\sum_{I\in\mathcal{A}}\left|I\right|_{m}\\
 & \le & (2N)^{-m}\epsilon\sum_{I\in D_{k}(\mathbb{R}^{m}),\, I\subset H_{1}}\left|I\right|_{m}=(2N)^{-m}\epsilon(2N)^{m}=\epsilon.
\end{eqnarray*}
Since this argument is valid for every $\epsilon>0$, we deduce from
Theorem \ref{thm:sandwich} that $\chi_{E}\in RI\left(\mathbb{R}^{m+1}\right)$
and $\int_{\mathbb{R}^{m+1}}\chi_{E}=0$, which completes the proof.
$\qed$

Here finally is the result about the integrability of functions which
are continuous everywhere except possibly on a very small set.
\begin{thm}
\label{thm:CtsAE}Suppose that the function $f:\mathbb{R}^{m}\to\mathbb{R}$
which is in $BB(\mathbb{R}^{m})$ and let $E$ be the set of points
$p\in\mathbb{R}^{m}$ at which $f$ is not continuous. If $E$ is
a very small subset of $\mathbb{R}^{m}$ then $f\in RI\left(\mathbb{R}^{m}\right)$. 
\end{thm}
\noindent \textit{Proof.} Let $R=\sup_{p\in\mathbb{R}^{m}}\left|f(p)\right|$.
Given an arbitrary positive number $\epsilon$, there exists some
$k\in\mathbb{N}$ and a finite collection $\mathcal{A}=\left\{ J_{1},J_{2},....,J_{N}\right\} $
of cubes in $D_{k}\left(\mathbb{R}^{m}\right)$ such that $E\subset\bigcup_{J\in\mathcal{A}}J$
and $\sum_{J\in\mathcal{A}}\left|J\right|_{m}<\epsilon/2R$. We define
the functions $g_{\epsilon}$ and $h_{\epsilon}$ by $g_{\epsilon}:=-R\sum_{J\in\mathcal{A}}\chi_{J}+\sum_{J\in D_{k}(\mathbb{R}^{m})\setminus\mathcal{A}}f\chi_{J}$
and $h_{\epsilon}:=R\sum_{J\in\mathcal{A}}\chi_{J}+\sum_{J\in D_{k}(\mathbb{R}^{m})\setminus\mathcal{A}}f\chi_{J}$.
(Of course only finitely many terms in each of these sums are non-zero.)
Since $f$ is continuous and bounded on each $J\in D_{k}(\mathbb{R}^{m})\setminus\mathcal{A}$,
we can use Theorem \ref{thm:CtySuffices} to obtain that $f\chi_{J}\in RI\left(\mathbb{R}^{m}\right)$
for each such $J$. Consequently $g_{\epsilon}$ and $h_{\epsilon}$
are both in $RI\left(\mathbb{R}^{m}\right)$. Furthermore, 
\[
\int_{\mathbb{R}^{m}}\left(h_{\epsilon}-g_{\epsilon}\right)=2R\int_{\mathbb{R}^{m}}\left(\sum_{J\in\mathcal{A}}\chi_{J}\right)=2R\sum_{J\in\mathcal{A}}\left|J\right|_{m}<\epsilon.
\]
The functions $g_{\epsilon}$ and $h_{\epsilon}$ also satisfy $g_{\epsilon}\le f\le h_{\epsilon}$.
Thus we can apply Theorem \ref{thm:sandwich} to obtain that $f\in RI(\mathbb{R}^{m})$.
$\qed$.

Here is a rather important example where we can use Theorem \ref{thm:CtsAE}. 
\begin{example}
\label{ex: YsimpleSet}Let $\left[a,b\right]$ be a bounded interval
and suppose that the functions $u:\left[a,b\right]\to\mathbb{R}$
and $v:\left[a,b\right]\to\mathbb{R}$ are both continuous on $\left[a,b\right]$
(this, as usual, meaning also that they are also one-sidedly continuous
at the endpoints). Suppose also that $u(x)\le v(x)$ for all $x\in[a,b]$.
Let $E=\left\{ (x,y):a\le x\le b:u(x)\le y\le v(x)\right\} $ and
let $f:\mathbb{R}^{2}\to\mathbb{R}$ be a bounded function which satisfies
$f(x,y)=0$ for all $(x,y)\in\mathbb{R}^{2}\setminus E$ and is continuous
at every interior point of $H$. Then $f$ is continuous at every
point of $(x,y)$ except possibly at points of the graphs $\left\{ (x,u(x)):a\le x\le b\right\} $
and $\left\{ (x,v(x)):a\le x\le b\right\} $ and the line segments
$\left\{ (a,y):u(a)\le y\le v(a)\right\} $ and $\left\{ (b,y):u(b)\le y\le v(b)\right\} $.
The first two of these four sets are very small subsets of $\mathbb{R}^{2}$
by Theorem \ref{thm:GraphIsVerySmall}. An obvious variant of this
theorem, with the roles of $x$ and $y$ interchanged shows that the
line segment sets are also very small subsets of $\mathbb{R}^{2}$.
Thus the union of all these sets is also very small and we can apply
Theorem \ref{thm:CtsAE} to obtain that $f\in RI(\mathbb{R}^{2})$. 

It is natural to ask whether the sufficient condition presented in
Theorem \ref{thm:CtsAE} for Riemann integrability of a function $f\in BB\left(\mathbb{R}^{m}\right)$
is also necessary. The answer to this question is no. However a certain
variant of this condition is necessary and sufficient. I am very tempted
to formulate that condition for you here, but for now I will resist
that temptation. It will arise rather naturally when you study the
a more sophisticated ``upgrade'' of the Riemann integral, namely
the Lebesgue integral. In this connection, see Exercise \ref{ex:MonFnIsRI}.

It is convenient to take note of the following rather obvious property
of very small subsets of $\mathbb{R}^{m}$.\end{example}
\begin{fact}
\label{fact:EqualAE}Suppose that $u\in RI\left(\mathbb{R}^{m}\right)$
and $f:\mathbb{R}^{m}\to\mathbb{R}$ is a bounded function such that
the set 
\[
E:=\left\{ p\in\mathbb{R}^{m}:u(p)\ne f(p)\right\} 
\]
 is a very small subset of $\mathbb{R}^{m}$. Then $f\in RI(\mathbb{R}^{m})$
and $\int_{\mathbb{R}^{m}}f=\int_{\mathbb{R}^{m}}u$. 
\end{fact}
\noindent \textit{Proof.} Let $R=\sup_{p\in\mathbb{R}^{m}}\left|f(p)-g(p)\right|$.
Then the fact that $f=u+\left(f-u\right)=u+\left(f-u\right)\chi_{E}$
implies that 
\[
u-R\chi_{E}\le f\le u+R\chi_{E}.
\]
We may set $g_{\epsilon}=u-R\chi_{E}$ and $h_{\epsilon}=u+R\chi_{E}$
for every $\epsilon>0$. Then $\int_{\mathbb{R}^{m}}\left(h_{\epsilon}-g_{\epsilon}\right)=2R\int_{\mathbb{R}^{m}}\chi_{E}=0$.
So Theorem \ref{thm:sandwich} and some obvious inequalities complete
the proof. $\qed$
\begin{rem}
\label{rem:IntegralFromAtoB}One immediate consequence of Fact \ref{fact:EqualAE},
Corollary \ref{cor:AllRectanglesMeasurable} and Theorem \ref{thm:CorollOfComposedFn}
is that for each $f\in RI(\mathbb{R})$ and $-\infty<a<b<\infty$
the four functions $f\chi_{(a,b)}$, $f\chi_{[a,b]}$, $f\chi_{[a,b)}$
and $f\chi_{(a,b]}$ are all in $RI(\mathbb{R})$ and their integrals
are all equal. ``Classically'' the same notation $\int_{a}^{b}f(x)dx$
is used to denote all four of these integrals. 
\end{rem}

We conclude this section by presenting two more very simple but often
useful properties of Riemann integrals. The first of these is a generalization
of Corollary \ref{cor:FinAditiv}.
\begin{fact}
\label{fact:Additivity}Suppose that the subsets $E_{1}$ and $E_{2}$
of $\mathbb{R}^{m}$ are disjoint or that $E_{1}\cap E_{2}$ is a
very small subset of $\mathbb{R}^{m}$, and $\chi_{E_{1}}$ and $\chi_{E_{2}}$
are both in $RI\left(\mathbb{R}^{m}\right)$. Then, for every $f\in RI\left(\mathbb{R}^{m}\right)$
the functions $f\chi_{E_{1}},$ $f\chi_{E_{2}}$ and $f\chi_{E_{1}\cup E_{2}}$
are in $RI\left(\mathbb{R}^{m}\right)$ and satisfy 
\begin{equation}
\int_{\mathbb{R}^{m}}f\chi_{E_{1}\cup E_{2}}=\int_{\mathbb{R}^{m}}f\chi_{E_{1}}+\int_{\mathbb{R}^{m}}f\chi_{E_{2}}.\label{eq:additive}
\end{equation}

\end{fact}

\noindent \textit{Proof.} Theorem \ref{thm:CorollOfComposedFn} ensures
that $f\chi_{E_{1}}$ and $f\chi_{E_{2}}$ are in $RI(\mathbb{R}^{m})$.
If $E_{1}\cap E_{2}=\emptyset$ then, as in Corollary  \ref{cor:FinAditiv},
we obtain the required conclusion from $\chi_{E_{1}\cup E_{2}}=\chi_{E_{1}}+\chi_{E_{2}}$.Otherwise
we need to use Fact \ref{fact:EqualAE} to ensure that, for $n=1,2$,
the function $f\chi_{E_{n}\setminus(E_{1}\cup E_{2})}$ is in $RI(\mathbb{R}^{m})$
and $\int_{\mathbb{R}^{m}}f\chi_{E_{n}\setminus\left(E_{1}\cap E_{2}\right)}=\int_{\mathbb{R}^{m}}f\chi_{E_{n}}$.
It is clear that $\chi_{E_{1}\setminus\left(E_{1}\cap E_{2}\right)}+\chi_{E_{2}\setminus\left(E_{1}\cap E_{2}\right)}=\chi_{\left(E_{1}\cup E_{2}\right)\setminus\left(E_{1}\cap E_{2}\right)}$.
So the function $g:=f\chi_{\left(E_{1}\cup E_{2}\right)\setminus\left(E_{1}\cap E_{2}\right)}$
is in $RI\left(\mathbb{R}^{m}\right)$ and $\int_{\mathbb{R}^{m}}g=\int_{\mathbb{R}^{m}}f\chi_{E_{1}}+\int_{\mathbb{R}^{m}}f\chi_{E_{2}}$.
But then yet another application of Fact \ref{fact:EqualAE} shows
that $f\chi_{E_{1}\cup E_{2}}$ is in $RI\left(\mathbb{R}^{m}\right)$
and $\int_{\mathbb{R}^{m}}f\chi_{E_{1}\cup E_{2}}=\int_{\mathbb{R}^{m}}g$,
which completes the proof. $\qed$

When $m=2$ the formula (\ref{eq:additive}) in more ``traditional''
notation, appears as 
\[
\iint_{E_{1}\cup E_{2}}f(x,y)dxdy=\iint_{E_{1}}f(x,y)dxdy+\iint_{E_{2}}f(x,y)dxdy.
\]
When $m=1$, if $-\infty<a<b<c<\infty$ and $f\in RI\left(\mathbb{R}\right)$
we can choose $E_{1}=\left(a,b\right)$ or $\left[a,b\right]$ and
$E_{2}=\left(b,c\right)$ or $\left[b,c\right]$ and in that case
(\ref{eq:additive}) can be rewritten as 
\begin{equation}
\int_{a}^{c}f(x)dx=\int_{a}^{b}f(x)dx+\int_{b}^{c}f(x)dx.\label{eq:intABC}
\end{equation}

\begin{fact}
\label{fact:AbsVal}Suppose that $f\in RI(\mathbb{R}^{m})$ and that
the subset $E$ of $\mathbb{R}^{m}$ satisfies $\chi_{E}\in RI(\mathbb{R}^{m})$.
Then $f\chi_{E}$ and $\left|f\right|\chi_{E}$ are in $RI\left(\mathbb{R}^{m}\right)$
and 
\begin{equation}
\left|\int_{\mathbb{R}^{m}}f\chi_{E}\right|\le\int_{\mathbb{R}^{m}}\left|f\right|\chi_{E}\le\sup_{p\in E}\left|f(p)\right|\int_{\mathbb{R}^{m}}\chi_{E}.\label{eq:eAbsVal}
\end{equation}

\end{fact}

\noindent \textit{Proof.} Theorem \ref{thm:CorollOfComposedFn} gives
us that $f\chi_{E}$ and $\left|f\right|\chi_{E}$ are in $RI\left(\mathbb{R}^{m}\right)$.
The pointwise inequalities $-\left|f\right|\chi_{E}\le f\chi_{E}\le\left|f\right|\chi_{E}$
and Theorem \ref{thm:Monotone} imply that $-\int_{\mathbb{R}^{m}}\left|f\right|\chi_{E}\le\int_{\mathbb{R}^{m}}f\chi_{E}\le\int_{\mathbb{R}^{m}}\left|f\right|\chi_{E}$
from which we obtain the first inequality in (\ref{eq:eAbsVal}).
The second inequality follows by applying Theorem \ref{thm:Monotone}
to the pointwise inequality $\left|f\right|\chi_{E}\le\sup_{p\in E}\left|f(p)\right|\chi_{E}$.
$\qed$. 

In the case where $f\in RI\left(\mathbb{R}^{m}\right)$ and $-\infty<a<b<\infty$
we can let $E$ be an open, closed or semiclosed interval whose endpoints
are $a$ and $b$ and use (\ref{eq:eAbsVal}) and then Corollary \ref{cor:AllRectanglesMeasurable}
to obtain that 
\begin{equation}
\left|\int_{a}^{b}f(x)dx\right|\le\int_{a}^{b}\left|f(x)\right|dx\le(b-a)\sup_{x\in\left(a,b\right)}\left|f(x)\right|.\label{eq:SimpEstm}
\end{equation}

\section{Methods for calculating Riemann integrals }

Given an explicit function $f\in RI\left(\mathbb{R}^{m}\right)$ we
can evaluate $\int_{\mathbb{R}^{m}}f$ at least approximately, by
calculating $L_{k}\left(f\right)$ or $U_{k}(f)$ for appropriate
(probably large) values of $k$. In some cases, certain properties
of $f$ may enable us to give an upper estimate for how close $U_{k}(f)$
or $L_{k}(f)$ is to $\int_{\mathbb{R}^{m}}f$ for given values of
$k$. Thus it may sometimes be possible, for each $\epsilon$, to
obtain a number $I_{\epsilon}$ which we can be sure satisfies $\left|I_{\epsilon}-\int_{\mathbb{R}^{m}}f\right|<\epsilon$.

In this section we wish to discuss other kinds of methods which sometimes
enable us to determine the value of $\int_{\mathbb{R}^{m}}f$ exactly.

\subsection{\label{sub:meo}The case $m=1,$ and the Leibniz-Newton formulae.}

\par \noindent Our starting point will be the following theorem. 
\begin{thm}
\label{thm:LeibNew1}Suppose that the function $g\in BB(\mathbb{R})$
is continuous at every point of the open bounded interval $\left(a,b\right)$.
Define the function $G:\left(a,b\right)\to\mathbb{R}$ by $G(t)=\int_{\mathbb{R}}g\chi_{(a,t)}$
for each $t\in\left(a,b\right)$. Then $G$ is differentiable and
satisfies 
\begin{equation}
G^{\prime}(t)=g(t)\mbox{ for every }t\in(a,b).\label{eq:PreLeibNewt}
\end{equation}
Furthermore 
\begin{equation}
\lim_{t\nearrow b}G(t)=\int g\chi_{(a,b)}\label{eq:GatB}
\end{equation}
 and 
\begin{equation}
\lim_{t\searrow a}G(a)=0.\label{eq:GatA}
\end{equation}

\end{thm}
\noindent \textit{Proof.}First we note that, for each $t$, the function
$g\chi_{\left(a,t\right)}\in RI\left(\mathbb{R}\right)$ by Theorem
\ref{thm:CtsAE} since it is continuous everywhere except possibly
at one or both of the points $a$ and $t$.

Fix an arbitrary point $t_{0}\in(a,b)$. For each real $h\ne0$ such
that $t_{0}+h\in\left(a,b\right)$ we have 
\[
G(t+h_{0})-G(t_{0})=\int_{\mathbb{R}}g\left(\chi_{(a,t_{0}+h)}-\chi_{(a,t_{0})}\right)=\alpha(h)\int_{\mathbb{R}}g\chi_{I(h)}
\]
where $\alpha(h)=1$ and $I(h)=[t_{0},t_{0}+h)$ when $h>0$ and,
if $h<0$, then $\alpha(h)=-1$ and $I(h)=[t_{0}+h,t_{0})$. So, for
each $h$ as above, 
\[
\frac{G(t_{0}+h)-G(t_{0})}{h}=\frac{\alpha(h)}{h}\int_{\mathbb{R}}g\chi_{I(h)}=\frac{1}{\left|h\right|}\int_{\mathbb{R}}g\chi_{I(h)}.
\]

Since $I(h)$ is an interval of length $\left|h\right|$ we have,
by Corollary \ref{cor:AllRectanglesMeasurable} and Remark \ref{rem:IntegralFromAtoB}
that $\int_{\mathbb{R}}g(t_{0})\chi_{I(h)}=g(t_{0})\left|h\right|$.
Therefore 
\begin{equation}
\frac{G(t_{0}+h)-G(t)}{h}-g(t_{0})=\frac{1}{\left|h\right|}\int_{\mathbb{R}}g\chi_{I(h)}-\frac{1}{\left|h\right|}\int_{\mathbb{R}}g(t_{0})\chi_{I(h)}=\frac{1}{\left|h\right|}\int_{\mathbb{R}}\left(g-g(t_{0})\right)\chi_{I(h)}.\label{eq:intf}
\end{equation}
Since $g$ is continuous at $t_{0}$, for each $\epsilon$ there exists
$\delta(\epsilon)$ such that $\left|g(x)-g(t_{0})\right|<\epsilon$
for every point $x\in(a,b)$ which satisfies $\left|x-t_{0}\right|<\delta(\epsilon)$.
If we choose $h$ such that $\left|h\right|<\delta(\epsilon)$, then
every $x\in I(h)$ satisfies $\left|x-t_{0}\right|<\left|h\right|<\delta(\epsilon)$.
It follows that $\left|g(x)-g(t_{0})\right|<\epsilon$ for all $x\in I(h)$
and so, using Fact \ref{fact:AbsVal} and then Corollary \ref{cor:AllRectanglesMeasurable}
and Remark \ref{rem:IntegralFromAtoB}, we have 
\[
\left|\int_{\mathbb{R}}(g-g(t_{0}))\chi_{I(h)}\right|\le\sup_{t\in I(h)}\left|g(t)-g(t_{0})\right|\int_{\mathbb{R}}\chi_{I(h)}\le\epsilon\left|h\right|
\]
and consequently $\left|\frac{1}{\left|h\right|}\int_{\mathbb{R}}\left(g-g(t_{0})\right)\chi_{I(h)}\right|\le\epsilon$.
Substituting this inequality in (\ref{eq:intf}) and applying this
argument for every $\epsilon>0$, we obtain that $\lim_{h\to0}\left(\frac{G(t_{0}+h)-G(t)}{h}-g(t_{0})\right)=0$
which is exactly what is needed to prove (\ref{eq:PreLeibNewt}). 

It remains to prove (\ref{eq:GatB}) and (\ref{eq:GatA}). Here again
we will use Fact \ref{fact:AbsVal} and then Corollary \ref{cor:AllRectanglesMeasurable}
and Remark \ref{rem:IntegralFromAtoB} to obtain that 
\[
\left|G(t)\right|\le\sup_{p\in\mathbb{R}}\left|g(p)\right|(t-a)
\]
and also to obtain that 
\[
\left|\int_{\mathbb{R}}g\chi_{(a,b)}-G(t)\right|=\left|\int_{\mathbb{R}}\left(g\chi_{\left(a,b\right)}-g\chi_{(a,t)}\right)\right|=\left|\int_{\mathbb{R}}g\chi_{[t.b)}\right|\le\sup_{p\in\mathbb{R}}\left|g(p)\right|(b-t)
\]
for all $t\in\left(a,b\right)$. These two inequalities immediately
give us (\ref{eq:GatA}) and (\ref{eq:GatB}) respectively. $\qed$
\begin{thm}
\label{thm:LeibNewt02}Suppose that the function $g\in BB(\mathbb{R})$
is continuous at every point of the bounded open interval $\left(a,b\right)$.
Suppose further that the function $F:\left[a,b\right]\to\mathbb{R}$
is continuous (included one-sidely at the endpoints) and is differentiable
and satisfies $F'(t)=g(t)$ at every point $t\in\left(a,b\right)$.
Then 
\begin{equation}
\int_{\mathbb{R}}g\chi_{(a,b)}=F(b)-F\left(a\right).\label{eq:LNF}
\end{equation}

\end{thm}

\noindent \textit{Proof.} Define the function $H:(a,b)\to\mathbb{R}$
by $H(t)=F(t)-G(t)$ where $G:\left(a,b\right)\to\mathbb{R}$ is the
function defined in the statement of Theorem \ref{thm:LeibNew1}.
Then $H'(t)=F'(t)-G'(t)=g(t)-g(t)=0$ for every $t\in(a,b)$. So,
by Lagrange's Theorem, $H$ equals a constant, say $c_{0}$ on $\left(a,b\right)$
and we can write $F(t)=G(t)+c_{0}$ for all $t\in(a,b)$. By the one-sided
continuity of $F$ at $a$ and (\ref{eq:GatA}) we have that $F(a)=c_{0}$.
Similarly (\ref{eq:GatB}) gives us that $F(b)=\int_{\mathbb{R}}g\chi_{(a,b)}+c_{0}$.
Subtracting these last two equations gives us (\ref{eq:LNF}). $\qed$

The formula (\ref{eq:LNF}) is often referred to as the Leibniz-Newton
formula or Newton-Leibniz formula. In the ``traditional'' presentation
of this result the integral on the left side of (\ref{eq:LNF}) is
denoted by $\int_{a}^{b}g(x)dx$ (cf.~also Remark \ref{rem:IntegralFromAtoB})
and $g$ is a function which is considered to be defined only on $\left(a,b\right)$
and perhaps at one or both of $a$ and $b$. Of course any values
taken by $g$ outside the interval $\left(a,b\right)$ are irrelevant
for Theorem \ref{thm:LeibNewt02}.

We refer to standard textbooks for a wide range of various techniques,
change of variables, special substitutions, integration by parts,
partial fractions, etc. etc. which, given a continuous bounded function
$g:\left(a,b\right)\to\mathbb{R}$ enable one to sometimes, but certainly
not always, find a function $F:\left[a,b\right]\to\mathbb{R}$ which
satisfies the hypotheses of \ref{thm:LeibNewt02}.

\subsection{\label{sub:SomePreps}Some preparations for the case $m\ge2$. Integrals
depending on a parameter. }

\par \noindent In this subsection we will present a special result
which will be needed later, when we seek to prove an important special
case of a general formula (in Corollary \ref{cor:SpecialFubini})
for calculating double integrals. 
\begin{thm}
\label{thm:DepPar}Suppose that $u:\left[a,b\right]\to\mathbb{R}$
and $v:\left[a,b\right]\to\mathbb{R}$ are continuous functions on
the closed bounded interval $\left[a,b\right]$ (one-sidedly at the
endpoints) which satisfy $u(x)\le v(x)$ for all $x\in[a,b]$. Let
$E=\left\{ (x,y):a\le x\le b,\, u(x)\le y\le v(x)\right\} $ and suppose
that the function $f:E\to\mathbb{R}$ is continuous on $E$, i.e.,
continuous in the sense that $\lim_{n\to\infty}f(x_{n},y_{n})=f\left(x_{*},y_{*}\right)$
for every point $\left(x_{*},y_{*}\right)\in E$ and every sequence
$\left\{ \left(x_{n},y_{n}\right)\right\} _{n\in\mathbb{N}}$ of points
in $E$ which converges to $\left(x_{*},y_{*}\right)$. Then the function
$\phi:\mathbb{R}\to\mathbb{R}$ defined by 
\[
\phi(x)=\left\{ \begin{array}{ccc}
{\displaystyle \int_{u(x)}^{v(x)}f(x,y)dy} & , & x\in[a,b]\\
\\
0 & , & x\in\mathbb{R}\setminus\left[a,b\right]
\end{array}\right.
\]
satisfies $\phi\in RI(\mathbb{R})$. 
\end{thm}

\begin{rem}
In the statement of this theorem we have reverted to the more ``traditional''
notation for the formula defining $\phi$. In this special case it
is more convenient than the usually preferable notation which we have
used more frequently in this document. According to that notation
we could also write $\phi(x)=\int_{\mathbb{R}}F_{x}$ where, for each
fixed $x\in\left[a,b\right]$, $F_{x}(y)=f(x,y)$ for each $y\in\left[u(x),v(x)\right]$,
and $F_{x}(y)=0$ for each $y\in\mathbb{R}\setminus[u(x),v(x)]$. 
\end{rem}

\begin{rem}
A somewhat more elaborate proof would show that the restriction of
$\phi$ to the interval $\left[a,b\right]$ is continuous, one-sidedly
at the endpoints. But we do not need that stronger result for our
purposes here. In the special case where $u$ and $v$ are constants,
this continuity of $\phi$ on $\left[a,b\right]$ follows almost immediately
from the uniform continuity of $f$ on the rectangle $E$.
\end{rem}
\noindent \textit{Proof of Theorem \ref{thm:DepPar}}.

Step 1: We start by considering the special case where $f(x,y)$ takes
a constant value, say $c$ on the set $E$. Then $\phi(x)=c\chi_{[a,b]}(v(x)-u(x))$
which is of course in $BB(\mathbb{R})$ and continuous except possibly
at the two points $a$ and $b$. So we are done.

Next we consider the case of a function which is not continuous on
$E$, namely a function of the form $f=c\chi_{I\cap E}$ where $I$
is a rectangle of the form $[\alpha,\beta)\times[\gamma,\delta)$.
In that case we have 
\[
\phi(x)=\chi_{[\alpha,\beta)\cap[a,b)}(x)\left(\min\left\{ v(x),\delta\right\} -\max\left\{ u(x),\gamma\right\} \right)
\]
which is again in $BB(\mathbb{R})$ and the only points where it can
fail to be continuous are the endpoints of the interval $[\alpha,\beta)\cap[a,b]$.
It follows that $\phi\in RI(\mathbb{R})$ in this case also.

Step 2: Now suppose that $u:\mathbb{R}^{2}\to\mathbb{R}$ is constant
on every $I\in D_{k}(\mathbb{R}^{2})$ for some fixed $k\in\mathbb{N}$
and $f=u\chi_{E}$. Then $f$ is a linear combination of finitely
many functions of the form $\chi_{I\cap E}$ and so, by the previous
step, $\phi$ is in $RI(\mathbb{R})$, since it is a linear combination
of functions in $RI(\mathbb{R})$.

Step 3: Since $E$ equipped with the usual euclidean metric on $\mathbb{R}^{2}$
is a compact metric space the function $f$ is uniformly continuous
on $E$. Thus, by considerations similar to those explained at the
beginning of the proof of Theorem \ref{thm:NeedUnifCty}, given any
$\epsilon>0$, there exists $k\in\mathbb{N}$ which is sufficiently
large to ensure that, for every $I\in D_{k}(\mathbb{R}^{m})$ which
intersects with $E$, every pair of points $p$ and $q$ in $I\cap E$
satisfies $\left|f(p)-f(q)\right|<\epsilon/M(b-a)$ where $M=\max\left\{ v(x)-u(x):x\in\left[a,b\right]\right\} $.
Thus we have 
\begin{equation}
\sup_{p\in I}f(p)\le\inf_{p\in I}f(p)+\frac{\epsilon}{M(b-a)}\label{eq:veq}
\end{equation}
for all such $I$. We define the functions $G_{\epsilon}:\mathbb{R}^{2}\to\mathbb{R}$
and $H_{\epsilon}:\mathbb{R}^{2}\to\mathbb{R}$ to be the finite sums
\[
G_{\epsilon}:=\sum_{I\in D_{k}(\mathbb{R}^{2}),\, I\cap E\ne\emptyset}\inf_{p\in I\cap E}f(p)\chi_{I\cap E}\mbox{ \,\,\,\ and\,\,\,\ }H_{\epsilon}:=\sum_{I\in D_{k}(\mathbb{R}^{2}),\, I\cap E\ne\emptyset}\sup_{p\in I\cap E}f(p)\chi_{I\cap E}.
\]
 It follows from (\ref{eq:veq}) that they satisfy 
\begin{equation}
G_{\epsilon}\le f\le H_{\epsilon}\le G_{\epsilon}+\frac{\epsilon}{M(b-a)}\mbox{ at every point of \ensuremath{E}. }\label{eq:2veq}
\end{equation}
Then we define $g_{\epsilon}:\mathbb{R}\to\mathbb{R}$ and $h_{\epsilon}:\mathbb{R}\to\mathbb{R}$
by 
\[
g_{\epsilon}(x)=\chi_{[a,b]}(x)\int_{u(x)}^{v(x)}G_{\epsilon}(x,y)dy\,\,\,\mbox{and}\,\,\, h_{\epsilon}(x)=\chi_{[a,b]}(x)\int_{u(x)}^{v(x)}H_{\epsilon}(x,y)dy.
\]
 Step 2 esures that the functions $g_{\epsilon}$ and $h_{\epsilon}$
are both in $RI(\mathbb{R})$. It follows from (\ref{eq:2veq}) (and
Theorem \ref{thm:Monotone}) that they also satisfy $g_{\epsilon}\le\phi\le h_{\epsilon}$
on $\mathbb{R}$ and 
\[
g_{\epsilon}(x)\le\phi(x)\le h_{\epsilon}(x)\le g_{\epsilon}(x)+\frac{\epsilon(v(x)-u(x))}{M(b-a)}\le\frac{\epsilon}{b-a}
\]
 for all $x\in\left[a,b\right]$. This in turn implies that $\int_{\mathbb{R}}(h_{\epsilon}-g_{\epsilon})=\int_{\mathbb{R}}\chi_{[a,b]}\left(h_{\epsilon}-g_{\epsilon}\right)\le\int_{\mathbb{R}}\frac{\epsilon\chi_{[a,b]}}{b-a}=\epsilon$.
Thus we have all conditions required to use Theorem \ref{thm:sandwich}
to obtain that $\phi\in RI(\mathbb{R})$. $\qed$

\subsection{Repeated integrals and the case $m\ge2$.}

\par \noindent In this subsection we shall see that, when a function
$f\in RI(\mathbb{R}^{m})$ satisfies certain additional hypotheses,
it is possible to reduce the calculation of the value of $\int_{\mathbb{R}^{m}}f$
to several (essentially $m$) calculations of integrals of functions
of one variable, which can, hopefully, be performed by the method
of Subsection \ref{sub:meo}, i.e., via Theorem \ref{thm:LeibNewt02}.

Perhaps when you read this subsection for the first time you should
only consider it in the special case where $m=2$. Indeed, for most
applications in this course it will suffice to consider the case $m=2$
and, occasionally, the case $m=3$. In those cases the following definition
will look a little simpler.
\begin{defn}
\label{def:RepeatedIntegral}For each integer $m\ge2$, let $RpI(\mathbb{R}^{m})$
denote the set of all functions $f\in BB\left(\mathbb{R}^{m}\right)$
which have the following properties:

(i) For each point $\vec{x}=\left(x_{1},x_{2},...,x_{m-1}\right)$,
the function $f_{\vec{x}}:\mathbb{R}\to\mathbb{R}$ defined by $f_{\vec{x}}(t)=f(x_{1},x_{2},...,x_{m-1},t)$
is an element of $RI\left(\mathbb{R}\right)$. 

\textit{(The point $\vec{x}=\left(x_{1},x_{2},...,x_{m-1}\right)$
is considered as a ``constant'' for this part of our definition,
and also during the integration of $f_{\vec{x}}$ which we will perform
in the next step.)}

(ii) The function $g_{f}:\mathbb{R}^{m-1}\to\mathbb{R}$ defined by
$g_{f}(x_{1},x_{2},...,x_{m-1})=\int_{\mathbb{R}}f_{\vec{x}}$ for
each $\vec{x}=\left(x_{1},x_{2},....,x_{m-1}\right)\in\mathbb{R}^{m-1}$
is an element of $RI(\mathbb{R}^{m-1})$.

\textit{(After performing this integration for each ``constant''
$\vec{x}$, we can now consider $\vec{x}$ as a ``variable'' point.)}

For each $f\in RpI\left(\mathbb{R}^{m}\right)$ we define the\textit{
repeated$\left(m-1,1\right)$ integral of $f$ on $\mathbb{R}^{m}$}
by

\begin{equation}
\mathcal{R}_{m-1,1}\int f=\int_{\mathbb{R}^{m-1}}g_{f}.\label{eq:dri}
\end{equation}
\end{defn}
\begin{rem}
\label{rem:ClassNotn}In more ``classical'' notation, the repeated
integral appearing in (\ref{def:RepeatedIntegral}) can be denoted
by 
\[
\int_{\mathbb{R}^{m-1}}\left(\int_{\mathbb{R}}f\left(x_{1},x_{2},...,x_{m-1},t\right)dt\right)dx_{1}dx_{2}...dx_{m-1}\,,
\]
where of course, in the framework of Riemann integration the integrals
over $\mathbb{R}$ and over $\mathbb{R}^{m-1}$ are in practice only
over bounded subsets of $\mathbb{R}$ and $\mathbb{R}^{m-1}$ since
the function $f$ vanishes on the compliment of a bounded subset of
$\mathbb{R}^{m}$. In particular, when $m=2$ we can rewrite the right
side of (\ref{eq:dri}) as $\int_{\mathbb{R}}\left(\int_{\mathbb{R}}f(x,y)dy\right)dx$.
If $m=3$ it can be written as $\int_{\mathbb{R}^{2}}\left(\int_{\mathbb{R}}f(x,y,z)dz\right)dxdy$.
\end{rem}

\begin{rem}
\label{rem:RepLin}By two successive applications of Theorem \ref{thm:rlinearity},
it is immediately evident that the set of functions $RpI\left(\mathbb{R}^{m}\right)$
is a vector space and that 
\begin{equation}
\mathcal{R}_{m-1,1}\int(\alpha f+\beta g)=\alpha\mathcal{R}_{m-1,1}\int f+\mathcal{\beta R}_{m-1,1}\int g\label{eq:rplin}
\end{equation}
for every $f$ and $g$ in $RpI\left(\mathbb{R}^{m}\right)$ and every
real $\alpha$ and $\beta$. Furthermore, for these functions $f$
and $g$, we can readily deduce, by two successive applications of
Theorem \ref{thm:Monotone}, that 
\begin{equation}
\mbox{If \,\,\ }f\le g\mbox{ \,\,\,\ then\,\,\,\ }\mathcal{R}_{m-1,1}\int f\le\mathcal{R}_{m-1,1}\int g\,.\label{eq:monrepint}
\end{equation}
\end{rem}
\begin{example}
\label{ex:RpI-PhiPsi}Suppose that $m=2$ and that $\phi$ and $\psi$
are functions in $RI(\mathbb{R})$. Then one can immediately verify
that the function $f:\mathbb{R}^{2}\to\mathbb{R}$ defined by $f(x,y)=\phi(x)\psi(y)$
for all $(x,y)\in\mathbb{R}^{2}$ is in $RpI(\mathbb{R}^{2})$ and
that 
\begin{equation}
\mathcal{R}_{1,1}\int f=\int_{\mathbb{R}}\phi\cdot\int_{\mathbb{R}}\psi.\label{eq:PhiPsi}
\end{equation}
In particular, if $\phi=\chi_{I_{1}}$ and $\psi=\chi_{I_{2}}$ where
$I_{1}$ and $I_{2}$ are both arbitrary intervals in $D_{k}\left(\mathbb{R}\right)$
for some $k$ then $f(x,y)=\chi_{I_{1}}(x)\chi_{I_{2}}(y)$ which
means that $f=\chi_{I_{1}\times I_{2}}$. So, using (\ref{eq:PhiPsi})
and also Theorem \ref{thm:cubeok} we see that each dyadic square
$J=I_{1}\times I_{2}$ in $D_{k}\left(\mathbb{R}^{2}\right)$ satisfies
\begin{equation}
\mathcal{R}_{1,1}\int\chi_{J}=\mathcal{R}_{1,1}\int\chi_{I_{1}\times I_{2}}=\int_{\mathbb{R}}\chi_{I_{1}}\int_{\mathbb{R}}\chi_{I_{2}}=2^{-2k}=\int_{\mathbb{R}^{2}}\chi_{J}\,.\label{eq:RepEqualsDouble}
\end{equation}
More generally, for any $m\ge2$, if $\phi\in RI(\mathbb{R}^{m-1})$
and $\psi\in RI(\mathbb{R})$ and $f:\mathbb{R}^{m}\to\mathbb{R}$
is defined by $f(x_{1},x_{2},....,x_{m})=\phi(x_{1},x_{2},...,x_{m-1})\psi(x_{m})$
for all $\left(x_{1},x_{2},...,x_{m}\right)\in\mathbb{R}^{m}$, then
$f\in RpI(\mathbb{R}^{2})$ and $\mathcal{R}_{m-1,1}\int f=\int_{\mathbb{R}^{m}}\phi\cdot\int_{\mathbb{R}}\psi$.

\medskip{}

\end{example}

\textcolor{blue}{}%

\bigskip{}

The following theorem gives rather general conditions which enable
double Riemann integrals to be calculated via repeated integration
of functions of one variable. All we need for the double integral
and repeated integral to be equal is that they should both simply
exist. But in fact this theorem is a quite special case of a much
more general result which also holds in the more general context of
Lebesgue integration (a difficult but important topic which is well
beyond the scope of these notes). That general result is known as
Fubini's Theorem. 
\begin{thm}
\label{thm:fubini1-1}Suppose that the function $f:\mathbb{R}^{2}\to\mathbb{R}$
is in both of the classes $RI\left(\mathbb{R}^{2}\right)$ and $RpI\left(\mathbb{R}^{2}\right)$.
Then 
\begin{equation}
\int_{\mathbb{R}^{2}}f=\mathcal{R}_{1,1}\int f.\label{eq:f11}
\end{equation}

\end{thm}

\begin{rem}
\label{rem:RepClassNot}In more ``classical'' notation, the formula
(\ref{eq:f11}), can be written as 
\[
\iint_{\mathbb{R}^{2}}f(x,y)dxdy=\int_{\mathbb{R}}\left(\int_{\mathbb{R}}f(x,y)dy\right)dx
\]
where again, in ``practice'', the integrations of $\mathbb{R}$
and over $\mathbb{R}^{2}$ are effectively over bounded subsets of
$\mathbb{R}$ and of $\mathbb{R}^{2}$. For example, if $f(x,y)$
vanishes outside a set $E$ of the form $E=\left\{ (x,y):a\le x\le b,\, u(x)\le y\le v(x)\right\} $
for real constants $a$ and $b$ and suitable functions $u$ and $v$
of one variable, then the ``classical'' way of writing (\ref{eq:f11})
is 
\[
\iint_{E}f(x,y)dxdy=\int_{x=a}^{b}\left(\int_{y=u(x)}^{v(x)}f(x,y)dy\right)dx.
\]
There are other alternative kinds of notation for the above repeated
integral (which I personally quite dislike). For example, in some
books you might encounter something like this: 
\[
\int_{a}^{b}dx\int_{u(x)}^{v(x)}f(x,y)dy.
\]

\textit{Proof of Theorem \ref{thm:fubini1-1}.} 

Step (i): For any $k\in\mathbb{N}$ suppose that $f=\chi_{J}$ where
$J$ is an arbitrary dyadic square in $D_{k}\left(\mathbb{R}^{2}\right)$.
Then, using Example \ref{ex:RpI-PhiPsi} and, in particular, (\ref{eq:RepEqualsDouble}),
we have that $f$ is in both $RI\left(\mathbb{R}^{2}\right)$ and
$RpI\left(\mathbb{R}^{2}\right)$ and satisfies the formula (\ref{eq:f11}). 

Step (ii): By the linearity properties of double integrals and of
repeated integrals (see Theorem \ref{thm:rlinearity} and Remark \ref{rem:RepLin})
we immediately deduce that every function $f$ which is a linear combination
$f=\sum_{n=1}^{N}\alpha_{n}\chi_{I_{n}}$ of characteristic functions
of dyadic squares of any size, i.e., a dyadic step function, also
satisfies $f\in RI\left(\mathbb{R}^{2}\right)\cap RpI\left(\mathbb{R}^{2}\right)$
and satisfies (\ref{eq:f11}). 

Step (iii): We are now ready to consider the general case of an arbitrary
function $f\in RI\left(\mathbb{R}^{2}\right)\cap RpI\left(\mathbb{R}^{2}\right)$. 

Since $f\in RI(\mathbb{R}^{2})$ we can apply Proposition \ref{prop:DyadicSandwich}
to obtain two sequences $\left\{ g_{k}\right\} _{k\in\mathbb{N}}$
and $\left\{ h_{k}\right\} _{k\in\mathbb{N}}$ of dyadic step functions
satisfying (\ref{eq:yhstzm}) and (\ref{eq:hbisl}). These properties,
together with Step (ii) and (\ref{eq:monrepint}), give us that 
\[
\int_{\mathbb{R}^{2}}g_{k}=\mathcal{R}_{1,1}\int g_{k}\le\mathcal{R}_{1,1}\int f\le\mathcal{R}_{1,1}\int h_{k}=\int_{\mathbb{R}^{2}}h_{k}\mbox{ for all }k\in\mathbb{N}.
\]
This means that the number $\mathcal{R}_{1,1}\int f$ has the property
mentioned in the last sentence of the statement of Proposition \ref{prop:DyadicSandwich}
and therefore must equal $\int_{\mathbb{R}^{2}}f$. $\qed$

In the light of the above proof, the generalization of this theorem
to larger values of $m$ should be a straightforward exercise:\end{rem}
\begin{xca}
\label{ex:Fubini}Prove that $\int_{\mathbb{R}^{m}}f=\mathcal{R}_{m-1,1}\int f$
whenever $f\in RI(\mathbb{R}^{m})\cap RpI(\mathbb{R}^{m})$.\end{xca}
\begin{rem}
Combining Theorem \ref{thm:fubini1-1} with an obvious modification
of it, where the roles of $x$ and $y$ are permuted, we can see that,
if $f\in RI(\mathbb{R}^{2})$ then it is possible to interchange the
order of integration in repeated integration, i.e., in ``classical''
notation, 
\begin{equation}
\int_{\mathbb{R}}\left(\int_{\mathbb{R}}f(x,y)dy\right)dx=\int_{\mathbb{R}}\left(\int_{\mathbb{R}}f(x,y)dx\right)dy.\label{eq:ChangOrdIntgn}
\end{equation}

However, when $f$ is not in $RI(\mathbb{R}^{2})$ it can happen that
both of the repeated integrals in (\ref{eq:ChangOrdIntgn}) exist,
but are not equal.
\end{rem}

Let us conclude this section by presenting a result which is an obvious
and immediate corollary of Theorem \ref{thm:fubini1-1}, Theorem \ref{thm:DepPar}
and Example \ref{ex: YsimpleSet}. It gives us a standard formula
which is often useful for calculating double integrals.
\begin{cor}
\label{cor:SpecialFubini}Suppose that $u:\left[a,b\right]\to\mathbb{R}$
and $v:\left[a,b\right]\to\mathbb{R}$ are continuous functions on
the closed bounded interval $\left[a,b\right]$ (one-sidedly at the
endpoints) which satisfy $u(x)\le v(x)$ for all $x\in[a,b]$. Let
$E=\left\{ (x,y):a\le x\le b,\, u(x)\le y\le v(x)\right\} $ and suppose
that the function $f:E\to\mathbb{R}$ is continuous on $E$. Then
$f\in RI\left(\mathbb{R}^{2}\right)\cap RpI\left(\mathbb{R}^{2}\right)$
and $\int_{\mathbb{R}^{2}}f=\mathcal{R}_{1,1}\int f$. In more ``classical''
notation, this last equation can be written as
\[
\iint_{\left\{ (x,y):a\le x\le b,u(x)\le y\le v(x,y)\right\} }f(x,y)dxdy=\int_{a}^{b}\left(\int_{y=u(x)}^{v(y)}f(x,y)dy\right)dx.
\]

\end{cor}

\begin{rem}
We should at least mention an analogue of this result for triple integrals,
namely the formula 
\[
\iiint_{E}f(x,y,z)dxdydz=\iint_{\Omega}\left(\int_{z=u(x,y)}^{v(x,y)}f(x,y,z)dz\right)dxdy
\]
 which holds, under suitable conditions, when $E$ is the set $\left\{ (x,y,z):(x,y)\in\Omega,\, u(x,y)\le z\le v(x,y)\right\} $.
Perhaps we will precisely formulate and prove it in a future version
of this document.
\end{rem}

\section{Some exercises}

Well so far we only have one exercise. There will be more:

\begin{xca}
\label{ex:MonFnIsRI}Suppose that $f:\mathbb{R}\to\mathbb{R}$ vanishes
outside some closed bounded interval $\left[a,b\right]$ and that
its restriction to $\left[a,b\right]$ is non-increasing. 

(i) Prove that $f\in RI(\mathbb{R})$ if $a=j_{1}2^{-k}$ and $b=j_{2}2^{-k}$
for integers $j_{1}$, $j_{2}$ and $k$ satisfying $j_{1}<j_{2}$. 

Do not try to use Theorem \ref{thm:CtsAE}. It does not apply to this
case. Hint: Prove and use the fact that, whenever $a\le\alpha<\beta<\gamma\le b$,
we have $\inf_{p\in[\alpha,\beta)}f(p)\ge\sup_{p\in[\beta,\gamma)}f(p)$. 

(ii) Use the result of part (i) and Theorem \ref{thm:sandwich} to
extend the result of part (i) to the case of all real numbers $a$
and $b$ satisfying $a<b$.

(iii)$^{*}$ This is a much more challenging and NOT compulsory exercise.
(It is at the level of the material of the course ``Real Functions''
which deals with Lebesgue measure and integration. Please do not spend
too much time on it, and certainly not at the expense of other exercises,
etc.) Justify my remark in part (i) by constructing a non-decreasing
function $f:\left[0,1\right]\to\mathbb{R}$ such that the subset $\Omega$
of $\left[0,1\right]$, consisting of those points where $f$ is not
continuous, satisfies $U_{k}\left(\chi_{\Omega},\mathbb{R}\right)=1$
for every $k\in\mathbb{N}$. (Obviously such an $f$ must be bounded.)
This exercise shows that, (as already mentioned above), the sufficient
condition for Riemann integrability expressed in Theorem \ref{thm:CtsAE}
is not a necessary condition. However, as also already mentioned above,
there is a modified version of that condition which \textbf{\uline{is}}
necessary and sufficient. It comes to light rather naturally when
one studies, ``davka'' (see \cite{R}), the Lebesgue integral. 
\end{xca}

\section{\label{sec:Equivalence}Our definition here is equivalent to the
``classical'' definition}

In this section will show that the definitions of Riemann integrals
$\int_{\mathbb{R}^{m}}f$ and the class $RI(\mathbb{R}^{m})$ of Riemann
integrable functions of $m$ variables, via dyadic intervals, squares
and cubes, which we have used in this document are equivalent to the
usual ``classical'' definitions appearing in most books and courses
about this topic. (See e.g., \cite{R} p.~105 (for $m=1)$ or \cite{S}
p.~48 (for general $m$) for these usual definitions.) As already
mentioned above, you can also see a rather shorter ``almost'' proof
of this fact for the case $m=1$ in the last section of \cite{ST}.

In our treatment here we will use the same formulation, terminology
and notation for these usual classical definitions as appears in \cite{S}
pp.~46--48, and we shall temporarily use the perhaps funny terminology
``\textbf{\textit{classically Riemann integrable function}}'' and
``\textbf{\textit{classical Riemann integral}}'' when it is necessary
to specify that we are working in terms of these usual definitions
rather than those of this document. The symbols $RI(\mathbb{R}^{m})$
and $\int_{\mathbb{R}^{m}}f$ will of course still have exactly the
same ``unclassical'' meanings as we gave them in the definitions
of Section \ref{sec:BasicProps}.

\subsection{A proof that semiclosed intervals/cubes can be replaced by closed
intervals/cubes in the definition of $RI(\mathbb{R}^{m})$}

As observed in Fact \ref{rem:ClassNotn} a given function $f\in BB(\mathbb{R}^{m})$
is in $RI(\mathbb{R}^{m})$ if and only if it satisfies the condition
(\ref{eq:iswz}). Our first task here will be to prove that this condition
is the same if, instead of taking suprema over the semiclosed intervals
$I$, we take the suprema over their closures $\overline{I}$. I.e.,
we have to show, for each $f\in BB(\mathbb{R}^{m})$ that 
\begin{equation}
\lim_{k\to\infty}\left(2^{-km}\sum_{I\in D_{k}(\mathbb{R}^{m})}\sup_{p,q\in\overline{I}}\left|f(p)-f(q)\right|\right)\le\lim_{k\to\infty}\left(2^{-km}\sum_{I\in D_{k}(\mathbb{R}^{m})}\sup_{p,q\in I}\left|f(p)-f(q)\right|\right).\label{eq:ponj}
\end{equation}

We will do this by showing that, for each $k_{0}\in\mathbb{N}$ and
each $\epsilon>0$, there exists $k_{*}$ such that, for each $n\ge k_{*}$,
\begin{equation}
2^{-nm}\sum_{J\in D_{n}(\mathbb{R}^{m})}\sup_{p,q\in\overline{J}}\left|f(p)-f(q)\right|\le2^{-k_{0}m}\sum_{I\in D_{k_{0}}(\mathbb{R}^{m})}\sup_{p,q\in I}\left|f(p)-f(q)\right|+\epsilon.\label{eq:hgq}
\end{equation}

Given arbitrary $f\in BB(\mathbb{R}^{m})$, $k_{0}\in\mathbb{N}$
and $\epsilon>0$, let $\mathcal{A}$ denote the, (necessarily finite)
collection of all $m$-cubes $I$ in $D_{k_{0}}(\mathbb{R}^{m})$
which have the property that $f$ does not vanish identically on $\overline{I}$
and, for each $n\ge k_{0}$, and each $I\in\mathcal{A}$ let $\mathcal{C}(I,n)$
denote the collection of all $2^{n-k_{0}}$ cubes $J$ in $D_{n}(\mathbb{R}^{m})$
which are contained in $I$. Then (\ref{eq:hgq}) is the same as 
\begin{equation}
2^{-nm}\sum_{I\in\mathcal{A}}\,\sum_{J\in\mathcal{C}(n,I)}\,\sup_{\,\, p,q\in\overline{J}}\left|f(p)-f(q)\right|\le2^{-k_{0}m}\sum_{I\in\mathcal{A}}\,\,\sup_{\, p,q\in I}\left|f(p)-f(q)\right|+\epsilon,\label{eq:BigHgq}
\end{equation}
since all terms which we have omitted in making the transition from
(\ref{eq:hgq}) to (\ref{eq:BigHgq}) are equal to $0$.

Let $N$ denote the number of $m$-cubes in $\mathcal{A}$ and let
$M=\sup_{p\in\mathbb{R}^{m}}\left|f(p)\right|$. We shall now apply
Proposition \ref{prop:ArbitraryRectangle} to each cube $I$ in $\mathcal{A}$
in order to obtain an integer $k_{*}$ such that, for each $n\ge k_{*}$
there exists an $m$-rectangle $A_{n}\left(I\right)$ with the following
properties:

(i) Its closure $\overline{A_{n}(I)}$ is contained in $I$, 

(ii) It satisfies $\left|I\right|_{m}-\left|A_{n}(I)\right|_{m}\le{\displaystyle \frac{\epsilon}{MN}}$.

(iii) It is the union of a finite collection, which we will denote
by $\mathcal{B}(I,n)$, of pairwise disjoint cubes in $D_{n}(\mathbb{R}^{m})$. 

Thus we have 
\begin{equation}
\left|I\right|_{m}-\sum_{J\in\mathcal{B}(I,n)}\left|J\right|_{m}<\frac{\epsilon}{MN}.\label{eq:hbt}
\end{equation}
Observe that the integer $k_{*}$ can be chosen to be same for all
choices of $I\in\mathcal{A}$. This is possible because there are
only finitely many cubes in $\mathcal{A}$. (But in fact any $k_{*}$
which has the required property for one cube $I\in\mathcal{A}$ will
also have it for every other $I'\in\mathcal{A}$, since the same operation
of translation which moves $I$ to $I'$ also moves $A_{n}\left(I\right)$
to $A_{n}(I')$.)

Property (i) of $A_{n}(I)$ ensures that each $J\in\mathcal{B}(I,n)$
satisfies $\overline{J}\subset I$ and we will use this fact in the
fourth line of the following calculation. 

We observe that $\left|I\right|_{m}=\sum_{J\in\mathcal{C}(I,n)}\left|J\right|_{m}$
and so, for those $I$ in the collection $\mathcal{A}$, the inequality
(\ref{eq:hbt}) can be rewritten as $\sum_{J\in\mathcal{C}(I,n)\setminus\mathcal{B}(I,n)}\left|J\right|_{m}<\epsilon/MN$.
We will use the equation and also the inequality of the preceding
sentence in the fifth line of the calculation for each $I\in\mathcal{A}$
which we now begin: 
\begin{eqnarray*}
 &  & 2^{-nm}\sum_{J\in\mathcal{C}(I,n)}\sup_{p,q\in\overline{J}}\left|f(p)-f(q)\right|=\\
 & = & \sum_{J\in\mathcal{C}(I,n)}\sup_{p,q\in\overline{J}}\left|f(p)-f(q)\right|\left|J\right|_{m}\\
 & = & \sum_{J\in\mathcal{B}(I,n)}\sup_{\, p,q\in\overline{J}}\left|f(p)-f(q)\right|\left|J\right|_{m}+\sum_{J\in\mathcal{C}(I,n)\setminus\mathcal{B}(I,n)}\sup_{p,q\in\overline{J}}\left|f(p)-f(q)\right|\left|J\right|_{m}\\
 & \le & \sum_{J\in\mathcal{B}(I,n)\,}\sup_{p,q\in I}\left|f(p)-f(q)\right|\left|J\right|_{m}+M\sum_{J\in\mathcal{C}(I,n)\setminus\mathcal{B}(I,n)}\left|J\right|_{m}\\
 & \le & \sup_{p,q\in I}\left|f(p)-f(q)\right|\left|I\right|_{m}+\frac{\epsilon}{N}=2^{-k_{0}m}\sup_{p,q\in I}\left|f(p)-f(q)\right|+\frac{\epsilon}{N}
\end{eqnarray*}
Thus we have shown that $2^{-nm}\sum_{J\in\mathcal{C}(I,n)}\sup_{p,q\in\overline{J}}\left|f(p)-f(q)\right|\le2^{-k_{0}m}\sup_{p,q\in I}\left|f(p)-f(q)\right|+\frac{\epsilon}{N}$
for each $I\in\mathcal{A}$. We now take the sum of these inequalities
as $I$ ranges over the $N$ elements of $\mathcal{A}$. This will
give us exactly (\ref{eq:BigHgq}) for each $n\ge k_{*}$ and therefore
completes the proof of (\ref{eq:ponj}).

We are now ready to show, in the following two subsections, that our
definitions of $RI(\mathbb{R}^{m})$ and $\int_{\mathbb{R}^{m}}f$
are equivalent to the ``classical'' definitions of Riemann integrable
functions on an $m$-rectangle and $m$-fold Riemann integrals.

\subsection{A proof that functions in $RI(\mathbb{R}^{m})$ are classically Riemann
integrable}

Suppose first that $f\in RI(\mathbb{R}^{m})$. Let $A$ be any closed
$m$-rectangle such that $f$ vanishes on $\mathbb{R}^{m}\setminus A$
and let $\phi$ be the restriction of $f$ to $A$. We shall show
that $\phi$ is classically Riemann integrable on $A$, and that the
classical Riemann integral $\int_{A}\phi(x^{1},x^{2},..,x^{m})dx^{1}dx^{2}...dx^{m}$
equals $\int_{\mathbb{R}^{m}}f$. 

According to Theorem 3-3 on p.\,48 of \cite{S}, in order to show
that $\phi$ is classically Riemannn integrable on $A$ it will suffice,
for each $\epsilon>0$, to find a partition $P$ of $A$ into non-overlapping
closed rectangles (i.e., $m$-rectangles) of the special form described
on p.~46 of \cite{S} such that 
\begin{equation}
U(\phi,P)-L(\phi,P)<\epsilon.\label{eq:ploomp}
\end{equation}
We shall fix some $k\in\mathbb{N}$ and take $P$ to consist of all
closed $m$-rectangles $A\cap\overline{J}$ as $J$ ranges over all
the cubes in $D_{k}(\mathbb{R}^{m})$ for which $A\cap\overline{J}$
is non-empty. This partition certainly has the special form described
on p.~46 of \cite{S}. The fact that $f\in RI(\mathbb{R}^{m})$ combined
with (\ref{eq:ponj}) implies that we can choose $k$ so that $2^{-km}\sum_{J\in D_{k}(\mathbb{R}^{m})}\sup_{p,q\in\overline{J}}\left|f(p)-f(q)\right|<\epsilon$,
which we will use in the third line of the following calculation which
will establish (\ref{eq:ploomp}). We will also use (\ref{eq:burz})
in its second line: 
\begin{eqnarray*}
 &  & U(\phi,P)-L(\phi,P)\\
 & = & \sum_{S\in P}\sup_{q\in S}\phi(q)v(S)-\sum_{S\in P}\inf_{q\in S}\phi(q)v(S)=\sum_{S\in P}\sup_{q,p\in S}\left|\phi(p)-\phi(q)\right|v(S)\\
 & = & \sum_{A\cap\overline{J}}\,\,\sup_{\in P\,\, p,q\in A\cap\overline{J}}\left|\phi(p)-\phi(q)\right|\left|A\cap\overline{J}\right|_{m}\le\sum_{J\in D_{k}(\mathbb{R}^{m})}\,\,\sup_{\,\, p,q\in\overline{J}}\left|f(p)-f(q)\right|\left|\overline{J}\right|_{m}<\epsilon.
\end{eqnarray*}
Now that we have established that $\phi$ is integrable on $A$ we
should of course also show that 
\begin{equation}
\int_{\mathbb{R}^{m}}f=\int_{A}\phi(x^{1},x^{2},..,x^{m})dx^{1}dx^{2}...dx^{m}.\label{eq:feg}
\end{equation}
This will follow from the following pointwise inequalities:
\begin{eqnarray*}
\sum_{A\cap\overline{J}\in P}\inf_{p\in A\cap\overline{J}}\phi(p)\chi_{J\cap A} & \le & \sum_{A\cap\overline{J}\in P}\inf_{p\in A\cap J}\phi(p)\chi_{J\cap A}\le f=f\chi_{A}\\
 & \le & \sum_{A\cap\overline{J}\in P}\sup_{p\in A\cap J}\phi(p)\chi_{J\cap A}\le\sum_{A\cap\overline{J}\in P}\sup_{p\in A\cap J}\phi(p)\chi_{J\cap A}.
\end{eqnarray*}
All the functions appearing in these inequalites are in $RI(\mathbb{R}^{m})$.
So we can apply Theorem \ref{thm:Monotone} and Corollary \ref{cor:AllRectanglesMeasurable}
to these inequalities to obtain that 
\begin{equation}
L(\phi,P)\le\int_{\mathbb{R}^{m}}f\le U(\phi,P).\label{eq:glump}
\end{equation}
We can now deduce (\ref{eq:feg}) from (\ref{eq:glump}), since (as
implied by the definition just before Theorem 3-3 on p.~48 of \cite{S})
\begin{equation}
L(\phi,P)\le\int_{A}\phi(x^{1},x^{2},..,x^{m})dx^{1}dx^{2}...dx^{m}\le U(\phi,P)\label{eq:jpt}
\end{equation}
for all partitions $P$ and we have been able to obtain a partition
$P$ satisfying (\ref{eq:ploomp}) for arbitrarily small values of
$\epsilon$. (This is similar to the reasoning used in the last statement
in the formulation of Theorem \ref{thm:sandwich}.) 

\bigskip{}

\subsection{A proof that classically Riemann integrable functions of $m$ variables
are in $RI(\mathbb{R}^{m})$}

It will be convenient to first introduce some notation. If $E=\prod_{j=1}^{m}[a_{j},b_{j}]$
is a closed $m$-rectangle, then we let $E^{\#}=\prod_{j=1}^{m}[a_{j},b_{j})$. 

By Corollary \ref{cor:AllRectanglesMeasurable} we know that 
\begin{equation}
\int_{\mathbb{R}^{m}}\chi_{E^{\#}}=\left|E^{\#}\right|_{m}=\left|E\right|_{m}=\int_{\mathbb{R}^{m}}\chi_{E}\label{eq:Esharp}
\end{equation}
 and so 
\[
\int_{\mathbb{R}^{m}}\chi_{E\setminus E^{\#}}=\int_{\mathbb{R}^{m}}\left(\chi_{E}-\chi_{E^{\#}}\right)=\left|E\right|_{m}-\left|E^{\#}\right|_{m}=0,
\]
 which shows that $E\setminus E^{\#}$ is a very small subset of $\mathbb{R}^{m}$.

Suppose that $A$ is a closed $m$-rectangle and that the function
$\phi:A\to\mathbb{R}$ is Riemann integrable on $A$. Let us define
the function $f:\mathbb{R}^{m}\to\mathbb{R}$ to equal $\phi$ on
$A$ and $0$ on $\mathbb{R}^{m}\setminus A$. We have to show that
$f\in RI(\mathbb{R}^{m})$ and that $\int_{\mathbb{R}^{m}}f=\int_{A}\phi(x^{1},x^{2},..,x^{m})dx^{1}dx^{2}...dx^{m}$. 

It will be convenient to introduce an auxiliary function $u:\mathbb{R}^{m}\to\mathbb{R}$
which equals $f$ on $A^{\#}$ and $0$ on $\mathbb{R}^{m}\setminus A^{\#}$.
Since $A\setminus A^{\#}$ is a very small subset of $\mathbb{R}^{m}$
and $f=u$ everywhere except possibly on $A\setminus A^{\#}$, it
will suffice (in view of Fact \ref{fact:EqualAE}) to show that $u\in RI(\mathbb{R}^{m})$
and that $\int_{\mathbb{R}^{m}}u$ equals the classical integral of
$\phi$ on $A$. 

Given $\epsilon>0$, there exists a partition $P$ which satisfies
(\ref{eq:ploomp}) and therefore 
\[
\sum_{S\in P}\left(\sup_{q\in S}\phi(q)-\inf_{q\in S}\phi(q)\right)\left|S\right|_{m}<\epsilon.
\]
Note that the $m$-rectangles $S^{\#}$ are pairwise disjoint and
their union is $A^{\#}$. These properties ensure that the functions
$g_{\epsilon}:=\sum_{S\in P}\inf_{q\in S^{\#}}u(q)\chi_{S^{\#}}$
and $h_{\epsilon}:=\sum_{S\in P}\sup_{q\in S^{\#}}u(q)\chi_{S^{\#}}$
satisfy $g_{\epsilon}\le u\le h_{\epsilon}$. Applying (\ref{eq:Esharp})
to each $S$, we obtain that 

\begin{eqnarray*}
\int_{\mathbb{R}^{m}}\left(h_{\epsilon}-g_{\epsilon}\right) & = & \sum_{S\in P}\left(\sup_{q\in S^{\#}}u(q)-\inf_{q\in S^{\#}}u(q)\right)\left|S^{\#}\right|_{m}\\
 & = & \sum_{S\in P}\left(\sup_{q\in S^{\#}}\phi(q)-\inf_{q\in S^{\#}}\phi(q)\right)\left|S\right|_{m}\\
 & \le & \sum_{S\in P}\left(\sup_{q\in S}\phi(q)-\inf_{q\in S}\phi(q)\right)\left|S\right|_{m}<\epsilon.
\end{eqnarray*}
We can therefore apply Theorem \ref{thm:sandwich} to deduce that
$u\in RI\left(\mathbb{R}^{m}\right)$. Finally, in order to show that
$\int_{\mathbb{R}^{m}}u$ has the required value, we observe that
\[
L(\phi,P)\le\int_{\mathbb{R}^{m}}g_{\epsilon}\le\int_{\mathbb{R}^{m}}u\le\int_{\mathbb{R}^{m}}h_{\epsilon}\le U(\phi,P)
\]
for our partition $P$ which also satisfies (\ref{eq:ploomp}). As
before, we can combine this with (\ref{eq:jpt}) and the fact that
$\epsilon$ can be arbitrarily small, to deduce that $\int_{\mathbb{R}^{m}}u=\int_{A}\phi(x^{1},x^{2},..,x^{m})dx^{1}dx^{2}...dx^{m}$.

\bigskip{}

Thank you very much for your attention, and for any thoughts about
this that you care to offer.

\textcolor{red}{}%

Michael Cwikel
\end{document}